\documentclass[12pt]{article}
\usepackage{xy}
\usepackage{amsfonts,amsmath}
\usepackage{amssymb,latexsym}
\usepackage[dvips]{epsfig}
\usepackage{amscd}
\usepackage{psfrag}

\usepackage[hmargin=3cm,vmargin=3.5cm]{geometry}

\title{Categorifications from planar diagrammatics}
\author{Mikhail Khovanov}
\date{August 30, 2010}
 
\newtheorem{prop}{Proposition}
\newtheorem{theorem}{Theorem}

\newcommand{\oplusop}[1]{{\mathop{\oplus}\limits_{#1}}}
 
\begin{document} 
     
\maketitle
\baselineskip 14pt
 

\def\R{\mathbb R}
\def\Q{\mathbb Q}
\def\Z{\mathbb Z}
\def\N{\mathbb N} 
\def\C{\mathbb C}
\def\l{\lbrace}
\def\r{\rbrace}
\def\o{\otimes}
\def\lra{\longrightarrow}
\def\RHom{\mathrm{RHom}}
\def\Id{\mathrm{Id}}
\def\mc{\mathcal}
\newcommand{\TenR}{\otimes}
\def\mf{\mathfrak} 
\def\Ext{\mathrm{Ext}}
\def\End{\mathrm{End}}
\newcommand{\oTen}[1]{\otimes_{#1}}
\def\gdim{{\mathrm{gdim}}}
\def\dmod{{\mathrm{-mod}}} 
\newcommand{\Hom}{{\rm Hom}}

\newcommand{\Ind}{{\rm Ind}}
\newcommand{\Res}{{\rm Res}}
      
\newcommand{\ig}[2]{\vcenter{\xy (0,0)*{\includegraphics[scale=#1]{#2}} \endxy}}
\newcommand{\igv}[2]{\vcenter{\xy (0,0)*{\reflectbox{\includegraphics[scale=#1, angle=180]{#2}}} \endxy}}
\newcommand{\igh}[2]{\vcenter{\xy (0,0)*{\reflectbox{\includegraphics[scale=#1]{#2}}} \endxy}}
\newcommand{\ighv}[2]{\vcenter{\xy (0,0)*{\includegraphics[scale=#1, angle=180]{#2}} \endxy}}
\newcommand{\ii}{\underline{\textbf{\textit{i}}}}
\newcommand{\igrotCW}[2]{\vcenter{\xy (0,0)*{\includegraphics[scale=#1, angle=90]{#2}} \endxy}}

\begin{abstract}  A diagrammatic presentation of functors and natural transformations  
and the virtues of biadjointness are discussed. 
We then review a graphical description of the category of Soergel bimodules 
and a diagrammatic categorification of positive halves of quantum groups. 
These notes are a write-up of Takagi lectures given by the author in Hokkaido 
University in June 2009. 
\end{abstract}


\psfrag{Xone}{if $i\cdot j=0$}  
\psfrag{Xtwo}{if $i\cdot j=-1$}
\psfrag{Xthree}{if $i \not= j$}
\psfrag{Xfour}{unless $i=k$ and $i\cdot j=-1$} 
\psfrag{A}{$\mathcal{A}$}
\psfrag{B}{$\mathcal{B}$}
\psfrag{X4}{$R(\nu)$}
\psfrag{X5}{$U^+$}
\psfrag{X6}{$\dim(U^+(\nu)), \ \nu\in \N^I$}

\section{Planarity of biadjointness} 

Adjoint functors, since their discovery by Daniel Kan~\cite{Kn} in 1958, have become 
quite ubiquitous in mathematics, with their universality well-documented already 
in the Wikipedia. We hope that biadjoint functors, which are pairs of 
functors $(E,F)$ such that $F$ is both right and left adjoint of $E$, will prove 
to be of importance as well.  

Let  us begin by reviewing the topological meaning of adjointness and biadjointness. 
We will depict a functor $F: \mc{A}\lra \mc{B}$ by a mark on a horizontal interval, 
with the half-intervals to the right and left of the mark labelled by $\mc{A}$ and $\mc{B}$, 
respectively. 

\begin{center}{\psfig{figure=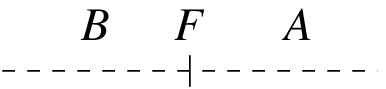,height=1cm}}\end{center}

Composition $F_n \dots F_1 : \mc{A}_1 \lra \mc{A}_{n+1}$ of functors 
$F_i: \mc{A}_i \lra \mc{A}_{i+1}$ is depicted by placing marks for $F_n, \dots, F_1$ 
in a row, with the intervals labelled by categories $\mc{A}_1, \dots, \mc{A}_{n+1}$ 
reading from right to left.  

\begin{center}{\psfig{figure=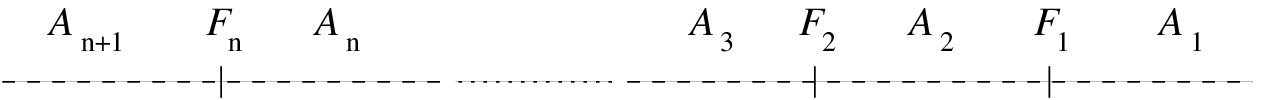,height=1.1cm}}\end{center}

The identity natural transformation $1_F : F\Rightarrow F$ is depicted by a 
vertical interval drawn in a rectangular region of the plane, with the two areas 
labelled by categories $\mc{A}$ and $\mc{B}$. 

\begin{center}{\psfig{figure=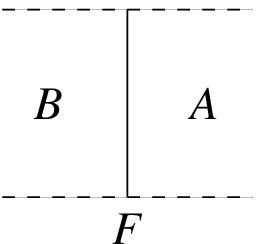,height=3cm}}\end{center}

A natural transformation $\alpha: F_1 \Rightarrow F_2$, where $F_1, F_2$ 
are functors from $\mc{A}$ to $\mc{B}$, is depicted by a dot on a vertical 
line. 

\begin{center}{\psfig{figure=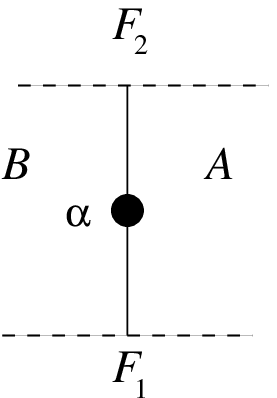,height=4cm}}\end{center}

A natural transformation $\alpha: F_n \dots F_1 \rightarrow G_m \dots G_1$ 
can be depicted by merging lines for the identity transformations of $F_n, \dots, F_1$ 
into a point and then splitting it into lines for the identities of $G_m, \dots, G_1$.  
 
\begin{center}{\psfig{figure=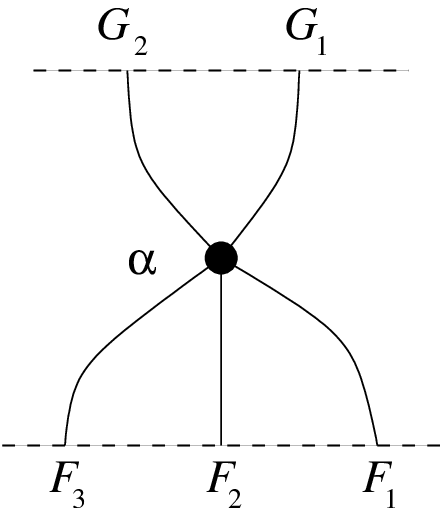,height=5cm}}\end{center}

It is useful not to label the identity functor $\Id_{\mc{A}}$ by anything and instead 
depict it by a horizontal interval labelled $\mc{A}$. Likewise, the identity 
natural transformation of this functor is denoted by a region labelled $\mc{A}$. 
With these rules, we can depict $\alpha: F_2 F_1 \Rightarrow \Id_{\mc{A}}$ by 
the following picture. 

\begin{center}{\psfig{figure=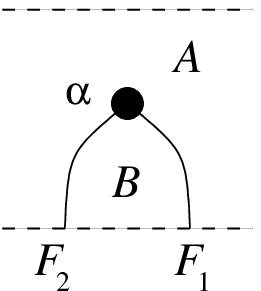,height=3cm}}\end{center}

\noindent 
Here $\mc{B}$ is the target category for functor $F_1$ and the source category 
for functor $F_2$. Similarly, below is a picture for $\alpha: \Id_{\mc{A}}\Rightarrow F$, 
where $F$ is an endofunctor of $\mc{A}$. 

\begin{center}{\psfig{figure=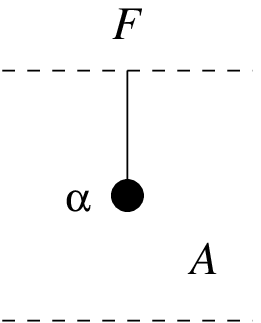,height=3cm}}\end{center}

The two possible types of compositions of natural transformations are depicted by either 
placing diagram in parallel or stacking them vertically; the next picture is an example of 
horizontal composition of three compatible natural transformations.  

\begin{center}{\psfig{figure=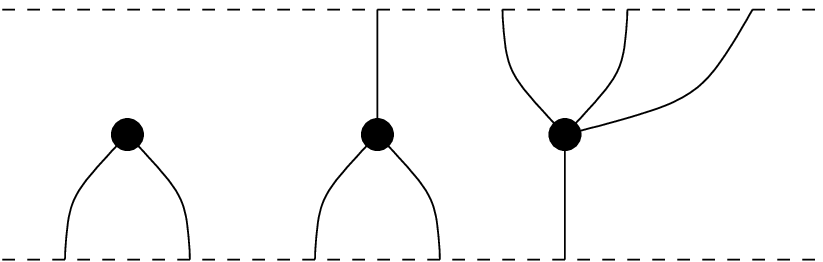,height=3cm}}\end{center}

An example of vertical composition of two natural transformations, with some identity 
transformations thrown in, is given below. 

\begin{center}{\psfig{figure=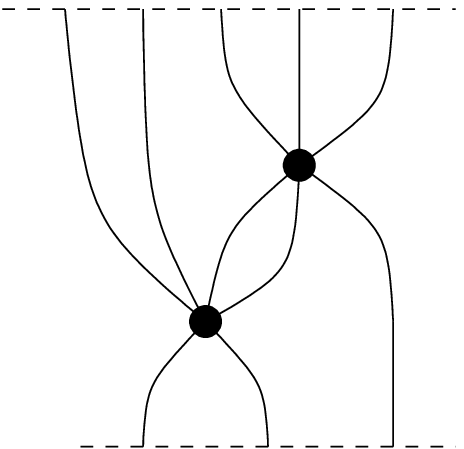,height=5cm}}\end{center}

It is important that strands always go ``up'', that is, strands are not yet 
allowed to have U-turns. Isotopies of strands are allowed, as long as they don't create U-turns. 
On some very informal level, these diagrams 
are analogous to planar projections of braids (braids don't have U-turns either). 

Thus, 2-dimensional planar pictures denote natural transformations, with 
regions of the picture labelled by categories, strands by functors and 
nodes by natural transformations. 

A planar diagram without boundary points (a closed diagram) determines an 
endomorphism of $\Id_{\mc{A}}$, where $\mc{A}$ is the label of the outside region, 
i.e., an element of the center of category $\mc{A}$. 

This setup, sometimes called string notation for 2-categories,  is Poincare dual to the more 
common one where categories are depicted by points, functors by arrows and 
natural transformations by 2-cells. The same setup can be used with any 
 2-category in lieu of the 2-category of natural transformations. 
Regions of the diagrams will be labelled by objects of the 2-category, 
edges by 1-morphisms, and nodes by 2-morphisms. 

A monoidal category is a 2-category with a single object. 
In this case we can avoid labelling regions, denote tensor product 
of objects of a monoidal category by a sequence of marks on a line and 
morphisms between tensor products by nodes with source and target 
arrows (there is a standard way of dealing with the case when the 
monoidal category is not strict; we omit the details). 
This notation is common in the diagrammatics for the 
representation categories of simple Lie algebras and quantum groups. 

We now come to adjointness. Functors  $F:\mc{A}\lra \mc{B}$ and 
$G: \mc{B}\lra \mc{A}$ are adjoint if there are isomorphisms 
$$ \Hom_{\mc{B}}(FX,Y) \cong \Hom_{\mc{A}}(X,GY)$$ 
functorial in $X\in \mc{A}$ and $Y\in \mc{B}$. An adjunction is equivalent to 
having the unit and the counit 
$$ \eta : \Id_{\mc{A}} \Rightarrow GF, \ \ \  \epsilon : FG \Rightarrow \Id_{\mc{B}}$$ 
natural transformations that satisfy the following equations: 
\begin{itemize} 
 \item  natural transformation $F \stackrel{1_F \circ \eta}{\lra} FGF 
\stackrel{\epsilon\circ 1_F}\lra F$ equals $1_F$, 
\item natural transformation $G \stackrel{\eta\circ 1_G}{\lra} GFG 
\stackrel{1_G\circ \epsilon}{\lra} G$ equals $1_G$. 
\end{itemize} 
Let us depict $\eta$ and $\epsilon$ by a cup and a cap diagram, respectively. 
To distinguish between functors $F$ and $G$ we equip the strands with up 
orientation near $F$ and down orientation near $G$. 

\begin{center}{\psfig{figure=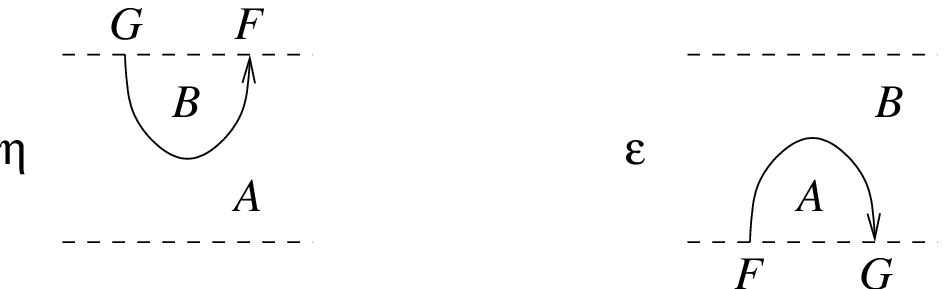,height=3.5cm}}\end{center}

The equations turn into relations on planar diagrams 

\begin{center}{\psfig{figure=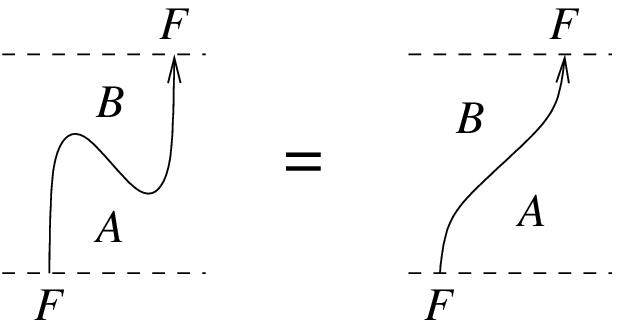,height=4cm}}\end{center}

\vspace{0.1in} 

\begin{center}{\psfig{figure=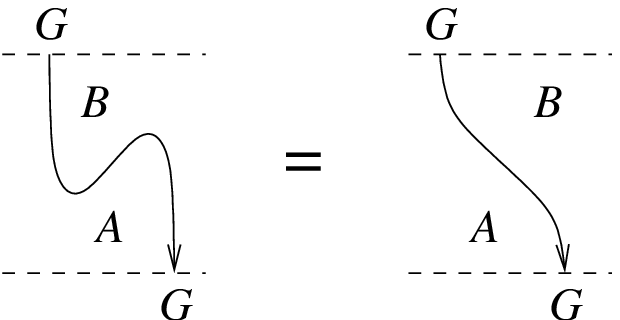,height=4cm}}\end{center}

These relations have a natural interpretation via isotopies of arcs, but give 
us only two out of four basic isotopy relations on oriented arcs in the plane. To get 
complete isotopy invariance, we assume that $F$ is also a right adjoint 
to $G$ and the natural transformations describing the second adjointness are 
fixed
$$ \underline{\eta} : \Id_{\mc{B}} \Rightarrow FG, \ \ \  \underline{\epsilon} : GF 
\Rightarrow \Id_{\mc{A}}.$$ 
They satisfy the relations   
\begin{itemize} 
\item natural transformation $F \stackrel{\underline{\eta}\circ 1_F}{\lra} FGF 
\stackrel{1_F\circ \underline{\epsilon}}{\lra} F$ equals $1_F$, 
 \item  natural transformation $G \stackrel{1_G \circ \underline{\eta}}{\lra} GFG 
\stackrel{\underline{\epsilon}\circ 1_G}\lra G$ equals $1_G$, 
\end{itemize} 
whose graphical interpretation is that of the two remaining types of arc isotopies:  

\begin{center}{\psfig{figure=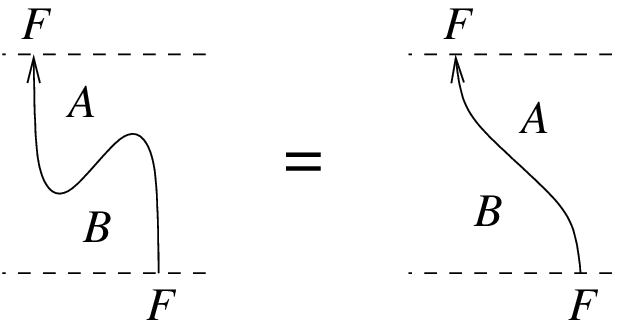,height=4cm}}\end{center}

\vspace{0.1in} 

\begin{center}{\psfig{figure=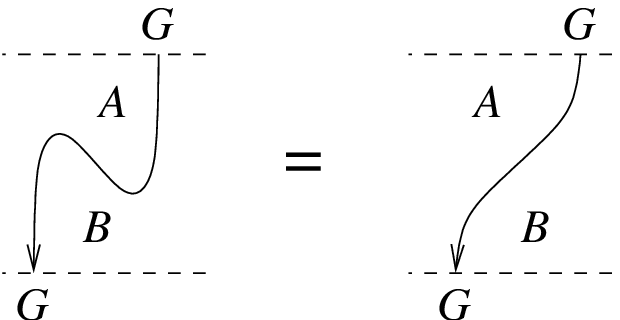,height=4cm}}\end{center}

Thus, given a biadjoint pair $(F,G)$ we can make sense out of any diagram 
built from  the four possible $(F,G)$-cups and caps. 
We say that $G$ is a biadjoint of $F$ and $F$ a biadjoint of $G$. A functor 
admitting a biadjoint is also called a \emph{Frobenius} functor, and a biadjoint 
pair $(F,G)$ is called a \emph{Frobenius pair}. Yet another terminology for a 
biadjoint pair is \emph{ambidextrous adjunction} \cite{Lau}. One of the  
first interesting examples of biadjoint pairs that appeared in mathematics 
(specifically, in the modular representation theory) 
were induction and restriction functors for inclusions of finite groups.
For a simpler example, take any invertible functor: its inverse is its biadjoint. 

One often works with many pairs of biadjoint functors, which requires compatibility, 
namely given biadjunctions $(F_1,G_1)$ and $(F_2,G_2)$ such that $F_2 F_1 $ 
is defined, the biadjunction for $(F_2F_1, G_1G_2)$ should be built in the 
obvious diagrammatic way via the composition rules. 

Biadjointness allows us to move dots on strands through cups and caps. 
Say, we have a transformation $\alpha: F_1\Rightarrow F_2$ and biadjoint 
pairs $(F_1,G_1), (F_2,G_2)$. There exists a unique transformation
$\alpha^{\ast}: G_2 \Rightarrow G_1$ that satisfies the equality

\begin{center}{\psfig{figure=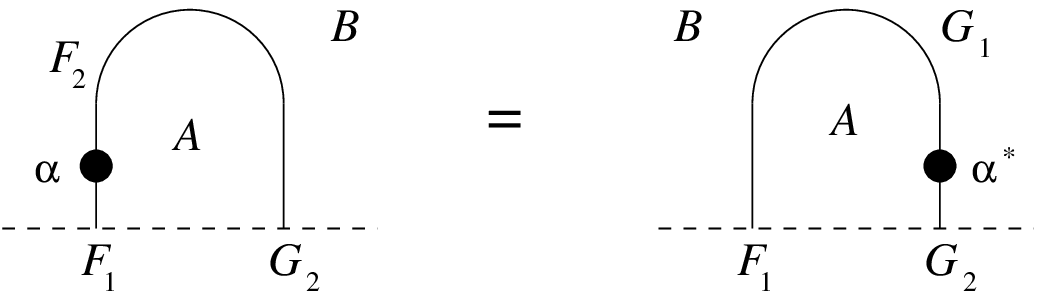,height=4cm}}\end{center}

${}^{\ast}\alpha: G_2 \Rightarrow G_1$ is defined similarly, via the equality 

\begin{center}{\psfig{figure=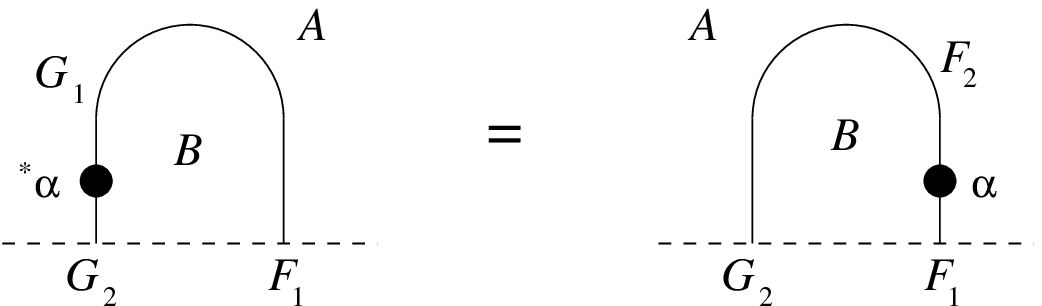,height=4cm}}\end{center}

From the topological viewpoint, it is convenient to require ${}^{\ast}\alpha=\alpha^{\ast}$
for all natural transformations $\alpha$ between our pairs of biadjoint functors, 
for then $\alpha$ always turns into $\alpha^{\ast}$ or $\alpha$ no matter what sequence of 
caps and cups it goes through. 
We refer to this property as \emph{cyclic biadjointness} (also called  \emph{even-handed structure} 
in~\cite{Bar}). 

With these assumptions, there is now a complete isotopy invariance for the string 
diagrams. Functors for these diagrams are selected from a collection of cyclic 
biadjoints that satisfy the above compatibility condition for their compositions. 
Any natural transformation between compositions of these functors can potentially 
be depicted. Some are glued out of the basic ones, others require introducing new 
nodes with multiple input and output strands. Once a new node for a natural transformation 
between composition appears, it may be used as a building block for more complicated 
planar diagrams. Due to cyclicity, these nodes can be isotoped in the plane. 

\emph{Remark:} If the cyclic order of functors around the node has a rotational symmetry, it 
might not extend to a symmetry of the node (but in some natural examples it does). 
For instance, given a node for $\alpha: F \Rightarrow F$ where $(F,F)$ is a biadjoint pair
(of course, $F: \mc{A}\lra \mc{A}$ is then an endofunctor), 
$\alpha^{\ast}: F\Rightarrow F$ does not have to be equal to $\alpha$.  

\vspace{0.1in} 

The above discussion generalizes, allowing us to depict elements of strict 2-categories 
with suitable duality properties, mirroring those of cyclic biadjointness. 

The planar interpretation of biadjointness has been a folklore for a number 
of years and appeared in~\cite{Muger, CKS}. String notation, in the case of 
2-categories with one object (monoidal categories), was introduced and made 
rigorous by Joyal and Street~\cite{JS}. Nowadays, it can be found on YouTube, in 
a series of videos "String diagrams" by TheCatsters.

\section{Biadjointness in topology and algebra}

Planar diagrams of lines labelled by functors and natural transformations can be 
thought of as suitably decorated 1-dimensional cobordisms embedded in $\R^2$. 
There is also a direct relation of biadjoint functors to cobordisms in all dimensions. 
Let $Cob_{n+2}$ be the 2-category of $(n+2)$-dimensional smooth cobordisms with 
boundary and corners. Closed $n$-dimensional manifolds $K$ are objects of $Cob_{n+2}$,  
while an $(n+1)$-dimensional cobordism $M$ is a 1-morphism from object $\partial_0 M$ 
to object $\partial_1 M$. An $(n+2)$-dimensional cobordism $N$ with corners is a 2-morphism 
from $\partial_0 N$ to $\partial_1 N$. Here $\partial_0N,\partial_1 N$ are $(n+1)$-manifolds 
with boundary, and the boundary of $N$ consists of 4 pieces: $\partial_0N,\partial_1 N$, 
and products $(\partial_0\partial_0N) \times [0,1]$,   $(\partial_1\partial_1N) \times [0,1].$
Corners of $N$ are four $n$-manifolds $\partial_i\partial_jN,$ $i,j\in\{0,1\}$. 

For each 1-morphism $M$ there is  the reverse 1-morphism $r(M)$ from $\partial_1 M$ to 
$\partial_0 M$ given by flipping $M$. We claim that $(M,r(M))$ constitutes a cyclic 
biadjoint pair. Multiply $M$ by an interval and then fan out 
the resulting $(n+2)$-manifold with corners so that it gives two $(n+2)$-cobordisms 
between $r(M) M$ and $\Id_{\partial_0M}$ (one in each direction) 
and two $(n+2)$-cobordisms 
between $M r(M)$ and $\Id_{\partial_1M}$. The relations on these four cobordisms 
are exactly the same as those satisfied by natural transformations 
$\eta, \epsilon, \underline{\eta}, \underline{\epsilon}$ of a biadjoint pair
(for more details and pictures see~\cite[Section 6.3]{FVIT}). 
In particular, given a 2-functor $F$ from $Cob_{n+2}$ into the 2-category of 
categories, functors, and natural transformations, the functors $F(M)$ and $F(r(M))$ 
are canonically biadjoint. 

If the target of $F$ is the 2-category of additive categories and additive functors, 
we say that $F$ is an $(n+2)$-dimensional TQFT with corners (or extended TQFT). 
Most of the time we require manifolds to be oriented, often additionally decorated 
by some structure, such as a spin structure, and sometimes place restrictions on 
topological type of manifolds and  cobordisms. There are interesting mathematically understood 
examples in dimensions $n+2=2, 3, 4.$ In dimension 2, they come from commutative 
Frobenius algebras. In dimension 3, the Witten-Reshetikhin-Turaev TQFT and its relatives 
emerge from Chern-Simons theory and representation theory of quantum groups, while  
the Rozansky-Witten TQFT comes from the derived category of coherent sheaves on 
holomorphic symplectic manifolds. Donaldson-Floer, Seiberg-Witten and Heegaard 
Floer theory (all closely interrelated) provide famous examples of 4D TQFTs. 

Any extended TQFT produces an abundance of biadjoint pairs, namely one pair $(F(M), F(r(M)))$ 
for each $(n+1)$-cobordism $M$. These biadjoint functors go between categories 
$F(K)$ assigned to closed $n$-manifolds $K$. When searching for TQFTs, it is useful 
to look for collections of categories that admit many biadjoint functors between each 
other. Examples of such collections include: 

\begin{itemize} 
\item Derived categories of coherent sheaves on a Calabi-Yau variety. 
Any sheaf on the product $X\times Y$ induces a pair of convolution 
functors between $D^b(X)$ and $D^b(Y)$. When $\dim(X)=\dim(Y)$, this 
is a biadjoint pair. When dimensions don't match, the functors are almost biadjoint 
(biadjoint up to a shift in the derived category). 
\item Categories of modules over finite-dimensional symmetric algebras and 
their derived categories. The functor of tensoring with a finitely-generated bimodule over 
symmetric algebras which is left and right projective has a bijadjoint. 
\item Fukaya-Floer categories of symplectic manifolds. 
\item Derived categories of sheaves on flag varieties. 
\item Various categories that appear in categorification of representations of Hecke 
algebras and quantum groups. For instance, functor $\mc{E}_i$ that categorifies 
generator $E_i$ of a simple Lie algebra/quantum group is biadjoint (or almost biadjoint) 
to the functor $\mc{F}_i$ categorifying generator $F_i$, see~\cite{CR, AL1, KL3, Rou2}. 
Some of the earliest examples of categorifications of Hecke algebra and $\mf{sl}(k)$ 
representations~\cite{GL, BFK, Su1, MS} 
came from highest weight categories of representations of $\mf{sl}(n)$, with 
$E_i$ and $F_i$ being translation functors (direct summands of tensor products with 
finite-dimensional modules) or Zuckerman functors. 
Biadjointness of translation functors was used already 
in~\cite{BG}. In Ariki's categorification~\cite{Ar} of irreducible $\mf{sl}(k)$ 
and affine $\mf{sl}(k)$ representations, functors $\mc{E}_i$ and $\mc{F}_i$ 
are biadjoint as well, being induction and restriction functors between 
Ariki-Koike-Cherednik cyclotomic quotients of affine Hecke algebras. 
\end{itemize} 

The issue of biadjointness in the Reshetikhin-Turaev-Witten TQFT is not prominent, since 
the category associated to the circle is semisimple $\C$-linear and functors associated 
to surfaces with boundary are $\C$-linear as well. Such functors are guaranteed to 
have biadjoints, which can be described in a simple combinatorial way. 
In contrast, the category associated to the circle in the Rozansky-Witten TQFT 
is a version of the derived category of coherent sheaves on a holomorphic symplectic 
manifold $X$, which is necessarily Calabi-Yau, so that one of the examples on the 
above list becomes relevant. Triangulated categories behind extended Donaldson-Floer and 
Ozsv\'ath-Szab\'o theories~\cite{LOT} are assigned to surfaces and  closely related to 
Fukaya-Floer categories of the representation variety of the fundamental group of a 
surface and of the symmetric power of a surface, respectively. Extended 4D TQFT 
that controls Ozsv\'ath-Szab\'o 3-manifold homology and 4-manifold invariants 
is being unravelled by Lipshitz, Ozsv\'ath and Thurston~\cite{LOT}.

\vspace{0.1in} 

We list a few nice features of biadjoint pairs: 

\begin{itemize} 
\item A biadjoint functor commutes with both limits and colimits. 
\item If a functor $F$ between additive categories has left or right adjoint, 
it is additive. In particular, in biadjoint pairs $(F,G)$ between additive categories 
both $F$ and $G$ are additive functors. 
\item In a biadjoint pair $(F,G)$ of functors between abelian categories both $F$ and $G$ 
are additive and exact, take projectives to projectives and injectives to injectives. 
\end{itemize} 

To an abelian category $\mc{A}$ there is assigned its Grothendieck group 
$G_0(\mc{A})$, an abelian group with generators given by symbols $[M]$ of objects of $\mc{A}$ 
and defining relations $[M]=[M']+[M'']$ for each short exact sequence 
$$ 0 \lra M''\lra M \lra M' \lra 0.$$ 
Often, it also makes sense to consider the group $K_0(\mc{A})$ generated by symbols 
of projective objects $[P]$ modulo relations $[P]=[P']+[P'']$ if $P\cong P'\oplus P''$. 
The obvious homomorphism $\phi_{\mc{A}}:K_0(\mc{A})\lra G_0(\mc{A})$ is, in general, 
neither injective nor surjective. The homomorphism takes $[P]$ to $[P]$, so this 
symbol notation might be ambiguous on projectives. 

A biadjoint pair $(F,G)$ induces homomorphisms 
$$[F] \  :  \  K_0(\mc{A}) \lra K_0(\mc{B}), \hspace{0.2in}
G_0(\mc{A}) \lra G_0(\mc{B}),  $$
$$[G] \  :  \ K_0(\mc{B}) \lra K_0(\mc{A}), \hspace{0.2in}  
G_0(\mc{B}) \lra G_0(\mc{A})$$
that commute with $\phi_{\mc{A}}, \phi_{\mc{B}}$. 
These homomorphisms take $[M]$ to $[FM]$, respectively $[GM]$. 
In contrast, an exact functor $F$ 
would induce a homomorphism $G_0(\mc{A}) \lra G_0(\mc{B})$ but not necessarily 
a homomorphism between $K_0$'s, since it might not take projectives to projectives. 
Likewise, a functor taking projectives to projectives would induce a homomorphism 
on $K_0$'s but not on $G_0$'s. Biadjoint pairs between abelian categories give the 
best behaving functors from this perspective. They can be thought of as categorifying 
pairs of adjoint operators. Namely, assume that $\mc{A},\mc{B}$ are $\Bbbk$-linear, 
over  a field $\Bbbk$, and hom spaces between objects are finite-dimensional over $\Bbbk$. 
This gives bilinear forms on $K_0(\mc{A}),$ $K_0(\mc{B})$ determined by 
$$ ([P],[Q])_{\mc{A}} := \dim \Hom_{\mc{A}}(P,Q) $$ 
for projectives $P, Q\in\mc{A}$, likewise for $\mc{B}$. Adjointness 
isomorphisms 
$$ \Hom_{\mc{B}}(FP,Q) \cong \Hom_{\mc{A}}(P,GQ),   \  \ \   
\Hom_{\mc{B}}(Q,FP) \cong \Hom_{\mc{A}}(GQ,P) $$
descend to relations 
$$([F]v,w)_{\mc{B}} = (v, [G]w)_{\mc{A}} , \  \   
    (w,[F]v)_{\mc{B}} = ([G]w,v)_{\mc{A}}, \  \  v\in K_0(\mc{A}), w\in K_0(\mc{B}),$$
showing that $[F]$ and $[G]$ become adjoint operators on real vector spaces 
$K_0(\mc{A})\otimes \R$ and  $K_0(\mc{B})\otimes \R$ relative to these two bilinear forms 
(in interesting examples the forms are often symmetric and nondegenerate).


\section{Diagrammatics for Soergel bimodules}

The Iwahori-Hecke algebra $H_n$ of the symmetric group $S_n$ has generators 
$T_i$, $1\le i\le n-1$, and defining relations
\begin{eqnarray*}
T_i^2 & = & (q-1)T_i + q, \\
T_iT_j & = & T_jT_i \ \  \mathrm{for}\ \  |i-j|\ge 2,  \\
T_iT_{i+1}T_i & = & T_{i+1}T_iT_{i+1},
\end{eqnarray*}
where $q$ is a formal parameter. 
In present-day mathematics it appears in two seemingly different ways: 

\vspace{0.06in} 

I) As a finite-dimensional quotient of the group algebra of the braid group, providing 
invariants of braid and links, when the latter are realized as closures of braids. 
Two braid closures produce the same link if they are related by a finite sequence of Markov moves. 
The algebraic counterpart of the closure operation is taking trace of an operator. 
The Ocneanu trace on the Hecke algebra behaves well under the Markov moves and 
can be normalized to produce an invariant of links known as the HOMFLY-PT polynomial 
(named with the initials of eight people who independently discovered it). More algebraically, 
the Hecke algebra is the commutant of the action of the quantum group $U_q(sl(k))$ 
on $V^{\otimes n}$ where $\dim (V)=k$ and $k\ge n$. 

\vspace{0.06in} 

II) As the endomorphism algebra of the representation of $GL(n,\mathbb{F}_q)$ 
induced from the trivial representation of the subgroup $\mc{B}$ 
of all invertible upper-triangular 
matrices. Here $GL(n,\mathbb{F}_q)$ is the finite group of invertible $n\times n$ 
matrices with coefficients in the finite field $\mathbb{F}_q$ with $q$ elements, hence
in this interpretation $q$ is no longer a formal variable. Case $q=1$ also fits into 
this framework, corresponding to $GL(n)$ over the one-element field, which is the 
symmetric group (and the Iwahori-Hecke algebra specializes to the group algebra of the symmetric 
group when $q=1$). 

\vspace{0.1in}

Construction I) of the Hecke algebra is fundamental for low-dimensional topology, 
construction II) and its generalizations is indispensable in representation theory. 
The second construction also leads to a categorification of the Iwahori-Hecke algebra. 
The functions on the finite set  $GL(n,\mathbb{F}_q)/\mc{B}$ become 
sheaves on the flag variety $GL(n)/\mc{B}$, either in the etale topology over 
a finite field, or  sheaves of $\C$-vector spaces on the flag variety over the field $\C$. 
Categorification of $H_n$ was constructed by Soergel~\cite{Soe1}, who also 
stated it in a sheaf-free language (which works for arbitrary Weyl groups). 
We describe it here in the symmetric group case.  

It is convenient to introduce $t=\sqrt{q}$ and view $H_n$ as a $\Z[t,t^{-1}]$-algebra. 
Then $H_n$ is also generated by $b_i=t^{-1}(T_i+1)$, $1\le i \le n-1$, 
with defining relations 
\begin{eqnarray}
b_i^2 & = & (t+t^{-1})b_i \label{eqn-bisq}\\
b_ib_j & = & b_jb_i \ \mathrm{for}\  |i-j|\ge 2 \label{eqn-bibj}\\
b_ib_{i+1}b_i + b_{i+1} & = & b_{i+1}b_ib_{i+1} + b_i \label{eqn-bibpbi}.
\end{eqnarray}

Soergel's categorification of $H_n$ is built out of bimodules over the polynomial 
ring $R=\C[x_1, \dots, x_n]$. The ring is graded, $\deg(x_i)=2$, and all the bimodules 
are graded as well. Grading shift $\{1\}$ up by one is an automorphism of the category of 
graded bimodules. For each $1\le i \le n-1$ let $R^i\subset R$ be the subring that 
consists of polynomials invariant under the transposition $x_i \leftrightarrow x_{i+1}$. Then 
$R \otimes_{R^i} R$ is a graded $R$-bimodule and we let 
$$B_i := R\otimes_{R^i} R\{-1\}.$$ 
Form the category $\mc{SC}_1$ whose objects are tensor products of $B_i$'s and 
morphisms are homomorphisms of bimodules. By adding finite  direct sums, 
grading shifts, restricting to grading-preserving homomorphisms, and 
forming the Karoubi envelope, one arrives at the category $\mc{SC}$. 
We call the objects of this category Soergel bimodules. 
One of Soergel's results is that $\mc{SC}$ categorifies the Hecke algebra $H_n$, i.e. 
the Grothendieck ring of $\mc{SC}$ is canonically isomorphic to $H_n$, 
$$ K_0(\mc{SC}) \cong H_n,$$ 
 with 
$[B_i]$, the symbol of $B_i$, going to $b_i$ under the isomorphism. 
Grading shift corresponds to multiplication 
by $t$: $[M\{1\}]= t [M]$. Multiplication in $H_n$ corresponds to the tensor product 
of bimodules,  
$$ [M]\cdot [N] \ \ := \ \ [M\otimes_R N].$$ 
Indecomposables in $\mc{SC}$ are in a bijection, up to grading shift, with 
elements $w$ of the symmetric group. Bimodule $B_e=R$, where $e$ is the trivial permutation, 
bimodules $B_{s_i}= B_i,$ where $s_i=(i,i+1)$, and, inductively, 
if $l(ww')=l(w)+l(w')$, $B_{ww'}$ is the only indecomposable summand of 
$B_w\otimes_R B_{w'}$ that is not isomorphic, up to a grading shift, to $B_u$ for some 
$u$ with $l(u)< l(w)+l(w')$.  Equivalently, 
$B_w$ is determined by the condition that it appears as a direct summand of
$B_{\ii}:=B_{i_1}\otimes B_{i_2}\otimes \dots \otimes B_{i_d}$, where $\ii =i_1\dots i_d$ and 
$s_{i_1}\dots s_{i_d}$ is a reduced presentation of $w$,
and does not appear as direct summand of any $B_{\ii}$, for sequences $\ii$ of length less than
$d=l(w)$.
   
Defining relations on $b_i$'s become isomorphisms of graded bimodules 
\begin{eqnarray}
B_i \TenR B_i & \cong & B_i\{ 1 \} \oplus B_i\{-1\},\\
B_i \TenR B_j & \cong & B_j \TenR B_i \ \mathrm{ for }\  |i-j|\ge 2,\\
(B_i \TenR B_{i+1} \TenR B_i) \oplus B_{i+1} & \cong & (B_{i+1} \TenR B_i \TenR B_{i+1}) 
\oplus B_i . \label{eqn-ipibimod}
\end{eqnarray}

The last isomorphism comes from decompositions 
\begin{eqnarray}
B_i \TenR B_{i+1} \TenR B_i \cong B_i \oplus (R \oTen{i,i+1} R\{-3\})\\
B_{i+1} \TenR B_i \TenR B_{i+1} \cong B_{i+1} \oplus (R \oTen{i,i+1} R\{-3\}).
\label{eqn-auxbimod2}
\end{eqnarray}
Here  $R \oTen{i,i+1} R$ is the tensor product of two $R$'s over the subring of 
$S_3$-invariants, with the action of $S_3$ by permutations of $x_i,x_{i+1},x_{i+2}$. 

When $n=3$, the Soergel category has 6 indecomposables, up to grading shifts. 
They are 
$$ R, \ \ B_1, \ \ B_2,  \ \ B_1\otimes_R B_2, \ \  B_2 \otimes_R B_1, \ \  
R \oTen{i,i+1} R\{-3\}.$$
Their images in the Grothendieck ring $K_0(\mc{SC})$ are 
$$ 1, \ \ b_1, \ \  b_2, \  \ b_1 b_2, \ \  b_2 b_1, \ \  b_1 b_2 b_1 - b_1 = b_2 b_1 b_2 -b_2.$$

N.~Libedinsky~\cite{Lib} presented $\mc{SC}$ in the $n=2$ case (and, more generally, 
in the so-called right-angled case, that we don't discuss) via generators and relations 
on morphisms. We now explain, following~\cite{EK}, a  
generalization of his presentation to an arbitrary $n$.  This approach provides a 
graphical description of homomorphisms between tensor products of Soergel bimodules.  

Tensoring with a Soergel bimodule is an endofunctor in the 
category of (graded) $R$-modules and homomorphisms between tensor products corresponds 
to natural transformations of functors. Start with the simplest Soergel bimodule $R$. 
The ring $R$ is commutative and multiplication 
by any element of $R$ is an endomorphism of the identity functor (tensoring with $R$), 
thus belongs to the center of the category of $R$-modules. 
We denote by a box labelled $i$ the multiplication by the generator $x_i$ of $R$. 

\begin{center}{\psfig{figure=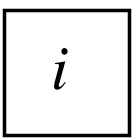,height=0.7cm}}\end{center}

One can think of this box as freely floating in a region. Boxes can float past each other 
in any direction (relative height change isotopy corresponds to commutativity of 
central elements). A collection of floating boxes denotes the product of corresponding 
generators. Any element of $R$ is a $\C$-linear combinations of products of boxes.  

A vertical line labelled $i$ will denote the identity endomorphism of the bimodule $B_i$. 
An important feature of $B_i$ is that it is selfadjoint (in the category of Soergel bimodules), 
i.e. left and right adjoint to itself, and this selfadjointness is cyclic. Therefore, we can 
introduce unoriented cup and cap diagrams labelled by $i$ to denote units and counits 
of biadjunctions (in pictures below, label $i$ is omitted). These diagrams have zero degree. 

Moreover, $B_i$ is a Frobenius object. Namely, there are homomorphisms 
$$R\lra B_i,  \ \ B_i \lra R, \ \ B_i \otimes B_i \lra B_i , \ \  B_i \lra B_i \otimes B_i $$ 
of degrees $1,1, -1, -1$ respectively, that we depict by 
\begin{equation} 
\ig{.3}{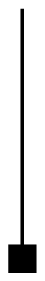} \hspace{0.4in} \igv{.3}{dot.eps} \hspace{0.4in} 
\igv{.3}{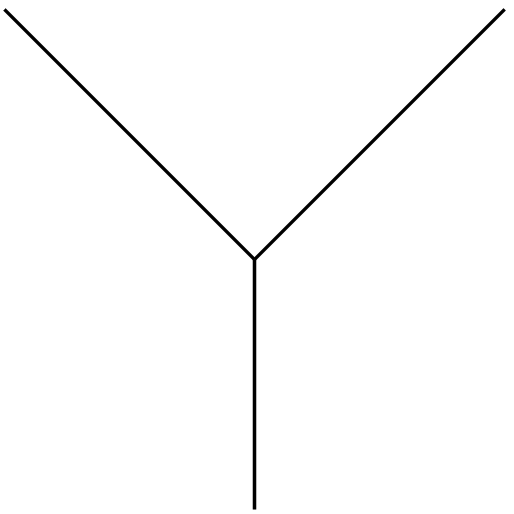}\hspace{0.4in} \ig{.3}{merge.eps}
\end{equation} 

and that satisfy the axioms of a Frobenius algebra object:  

\begin{equation}\ig{.2}{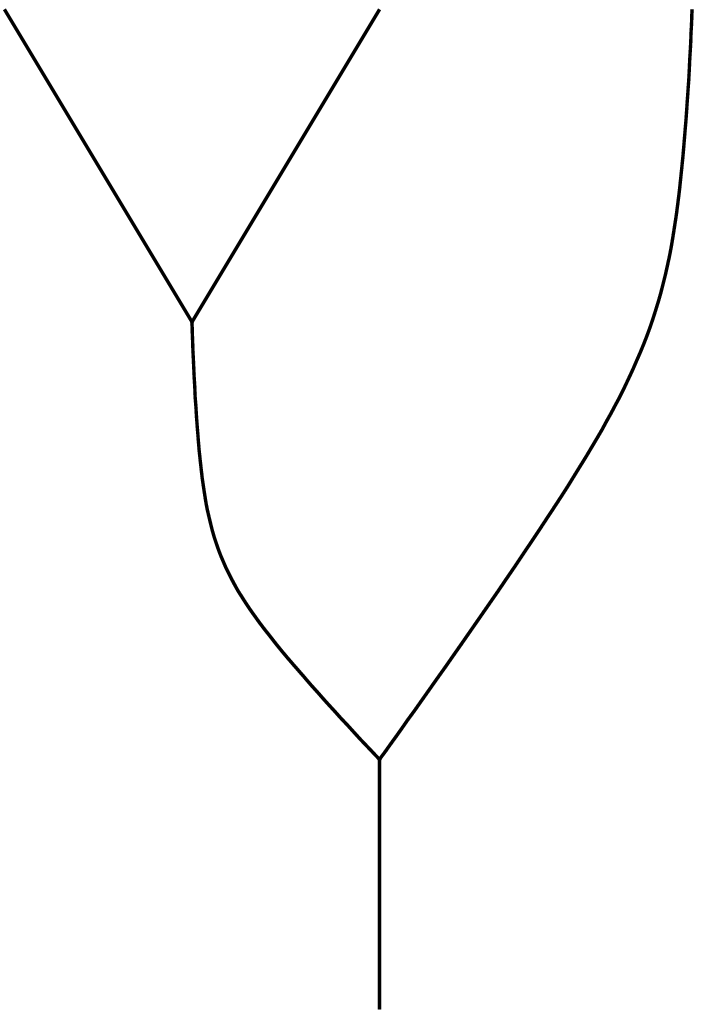} = \igh{.2}{doubleMerge.eps}  \hspace{0.4in} 
\igv{.2}{doubleMerge.eps} = \ighv{.2}{doubleMerge.eps} 
\end{equation}
\vspace{0.1in} 
\begin{equation}\ig{.2}{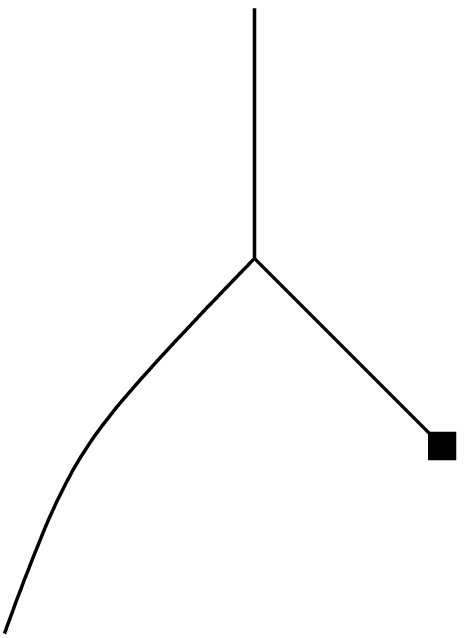} = \ig{.5}{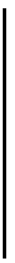} = \igh{.2}{splitDot.eps}
\hspace{0.4in} \igv{.2}{splitDot.eps} = \ig{.5}{line.eps} = \ighv{.2}{splitDot.eps}
\end{equation}
\vspace{0.1in} 
\begin{equation} \ig{.25}{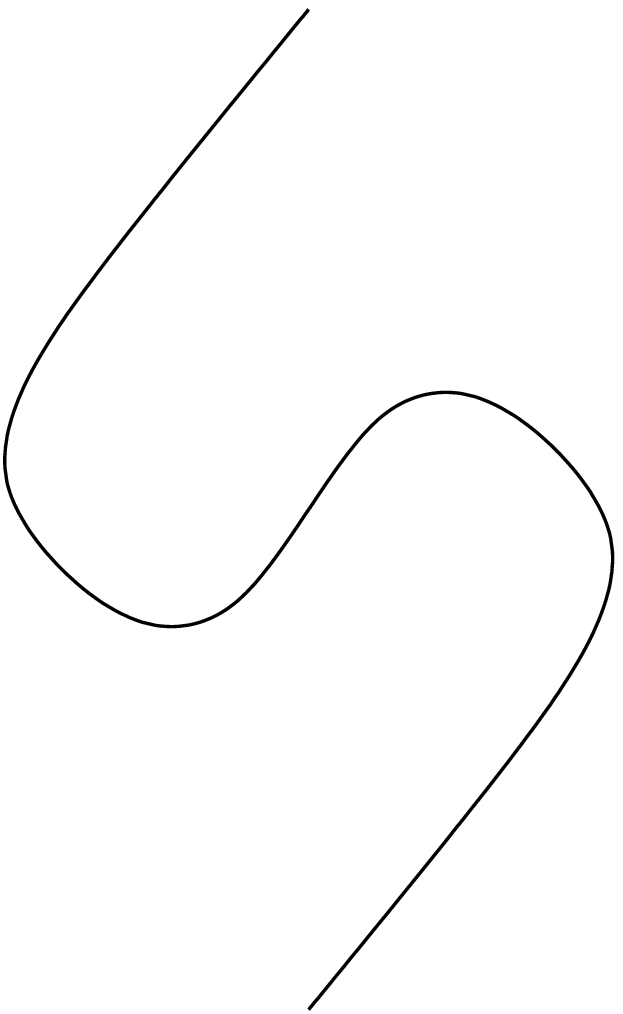} \;\; = \;\; \ig{.2}{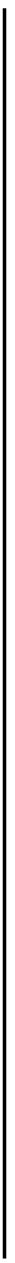} \;\; = \;\; 
 \igh{.25}{cupCap.eps}\hspace{0.3in} 
\ig{.3}{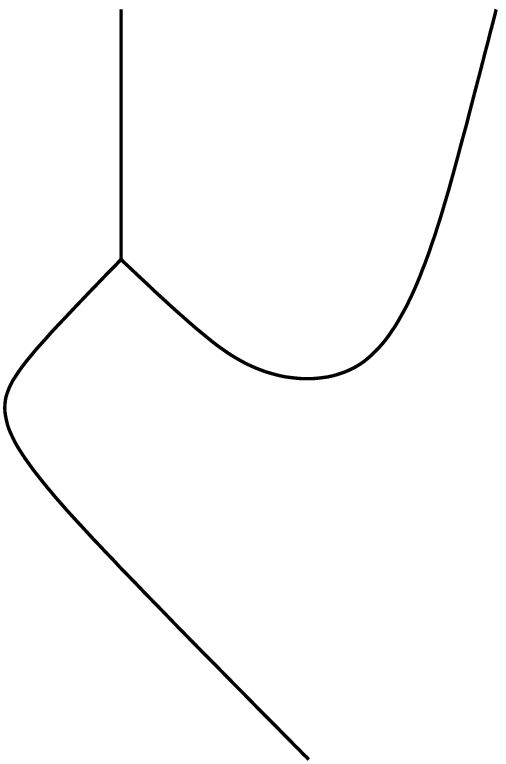} \  = \ \ig{.3}{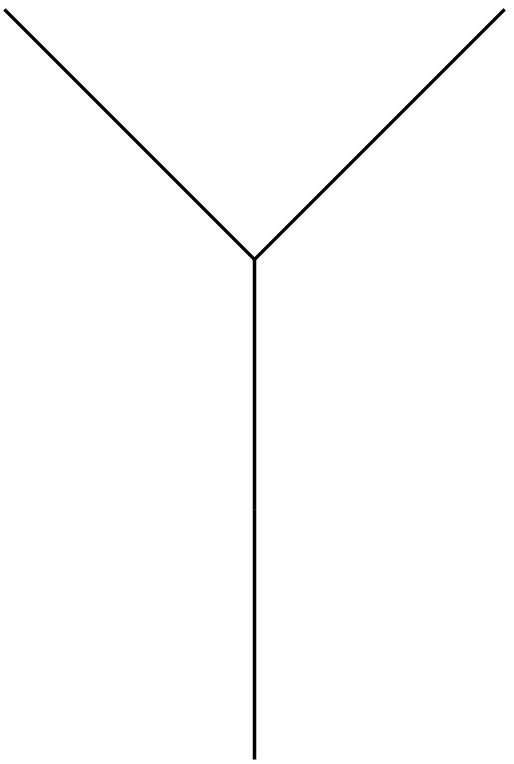} \  = \ \igh{.3}{splitCup.eps}
\end{equation}
\vspace{0.1in} 
\begin{equation}\ig{.3}{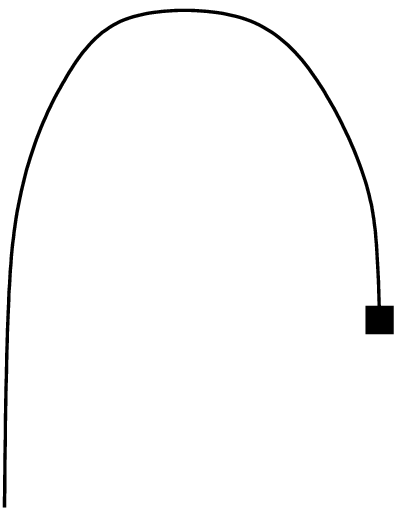} \  =  \ \igv{.4}{dot.eps}  \ =  \ \igh{.3}{capDot.eps}
\label{twistDot1}\end{equation}
\vspace{0.1in} 
\begin{equation}    \igv{.3}{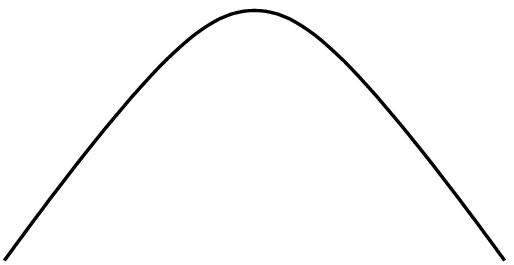} \ \ = \ \  \igv{.3}{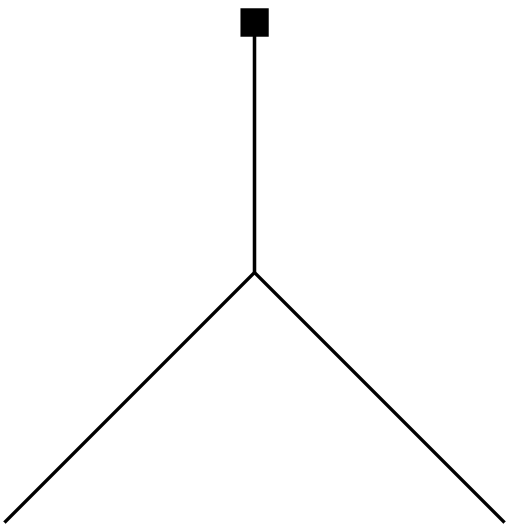} \hspace{0.2in} 
\ig{.4}{cap.eps} \ \ = \ \  \ig{.4}{dotSplit.eps}   \label{twistDot4}\end{equation}

We list other relations that involve strands of one color only:

\begin{eqnarray}
\ig{.4}{multxi.eps}\ \ig{.6}{line.eps}\ +\ \ig{.4}{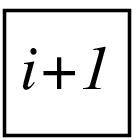}\ \ig{.6}{line.eps}\ &\ =
  \ &\ \ig{.6}{line.eps}\ \ig{.4}{multxi.eps}\ +\ \ig{.6}{line.eps}\ \ig{.4}{multxp.eps}
 \label{eq-slide1} \\    \ \  &  &  \ \  \nonumber  \\
\ig{.4}{multxi.eps}\ \ig{.4}{multxp.eps}\ \ig{.6}{line.eps}\ &\ =\ &\ \ig{.6}{line.eps}\ 
  \ig{.4}{multxi.eps}\ \ig{.4}{multxp.eps} \label{eq-slide2} \\
  \ \  &  &  \ \   \nonumber \\  
\ig{.4}{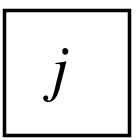}\ \ig{.6}{line.eps}\ &\ =\ &\ \ig{.6}{line.eps}\ \ig{.4}{multxj.eps}, 
\hspace{0.2in} |i-j|>1. 
\label{eq-slide3}
\end{eqnarray}

These come directly from the definition of $B_i$. Furthermore, we impose 

\begin{equation}\ig{.5}{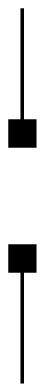} \ = \  \ig{.4}{multxi.eps}\ \ig{.75}{line.eps} \ -
   \  \ig{.75}{line.eps}\ \ig{.4}{multxp.eps} \ \ \ \ \ \  \ \ \ \ \ \ \ \ \ 
\ig{.5}{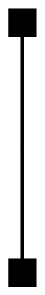} \ = \ \ig{.4}{multxi.eps} \ - \ \ig{.4}{multxp.eps} 
 \end{equation}

\begin{equation}\ig{.35}{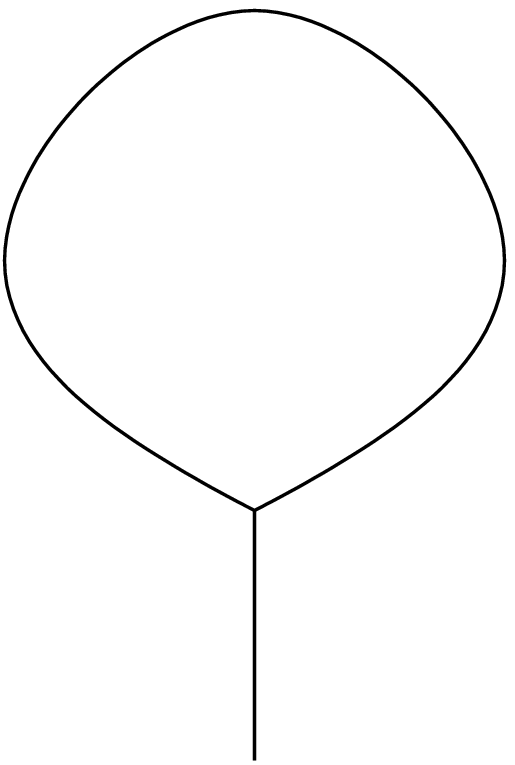} \ \ = \ \ 0 
\end{equation}

These relations imply, in particular, 

\begin{align} \ig{.4}{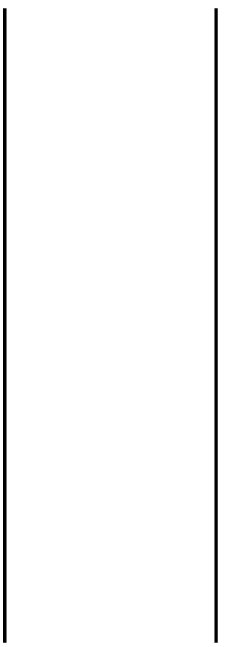}& \ \  = \ \  \ig{.4}{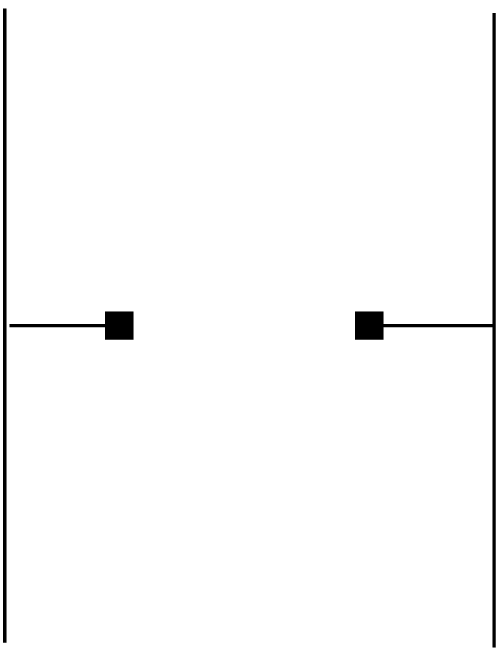} \ \  = \ 
\ig{.4}{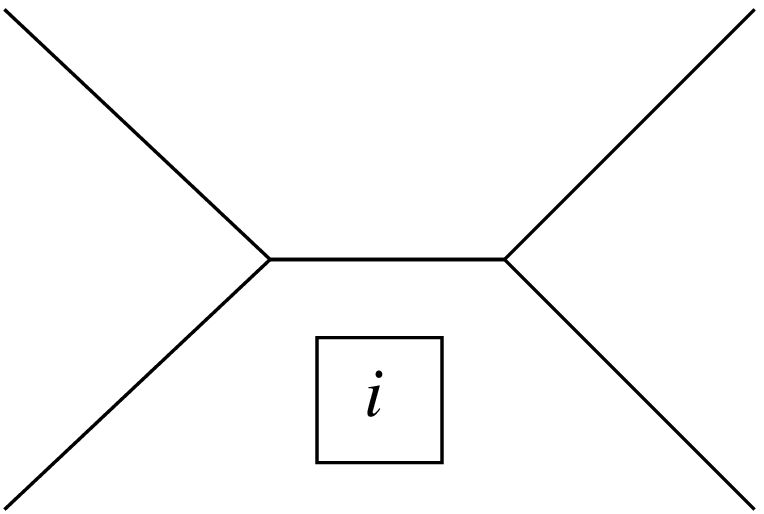} \ - \ \ig{.4}{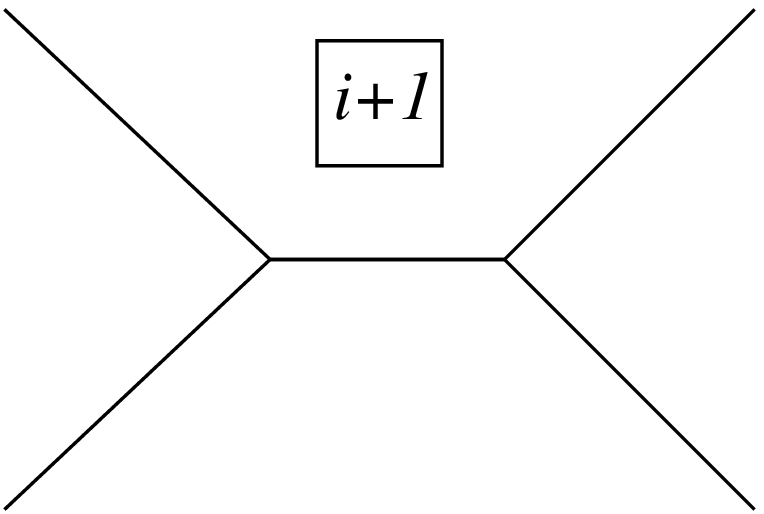} \label{twoLines}\\   \ \nonumber   \\
 &= \ \ig{.4}{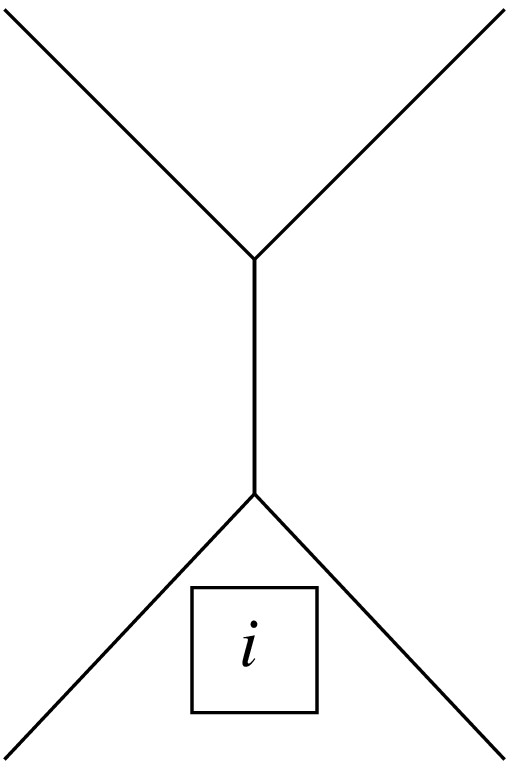} \ -\ \ig{.4}{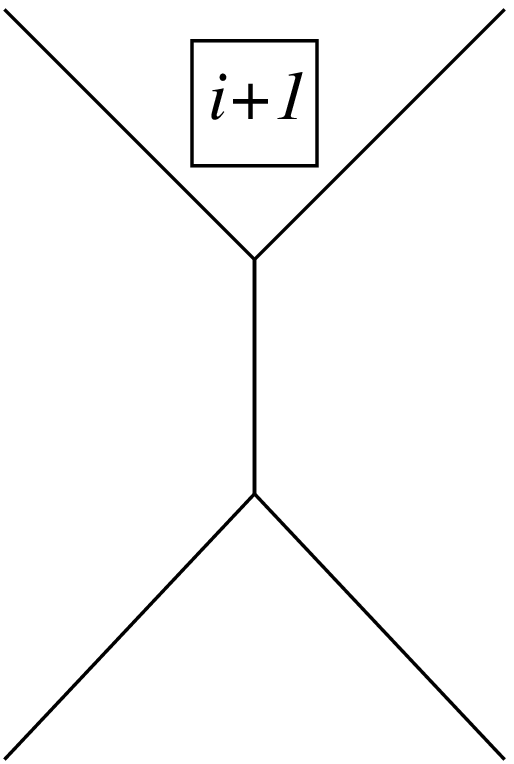} \nonumber \end{align}

leading to the decomposition $B_i\otimes B_i \cong B_i\{1\}\oplus B_i\{-1\}$. 

\vspace{0.2in} 

We now describe generators and relations for interactions of two adjacent colors $i$ and $i+1$. 
Below, thin lines represent $B_i$ and thick lines $B_{i+1}$. 
The six-valent vertex denotes the composition 
 $$ B_i\otimes B_{i+1}\otimes B_i \lra R_{\otimes_{i,i+1}} R\{-3\} \lra 
   B_{i+1}\otimes B_i \otimes B_{i+1}$$ 
of projection from the tensor product of 3 bimodules onto its most interesting direct summand 
and inclusion of this summand into the other triple tensor product. It turns out that the 
6-valent vertex possesses full rotational invariance: 

\begin{equation} 
\ig{.32}{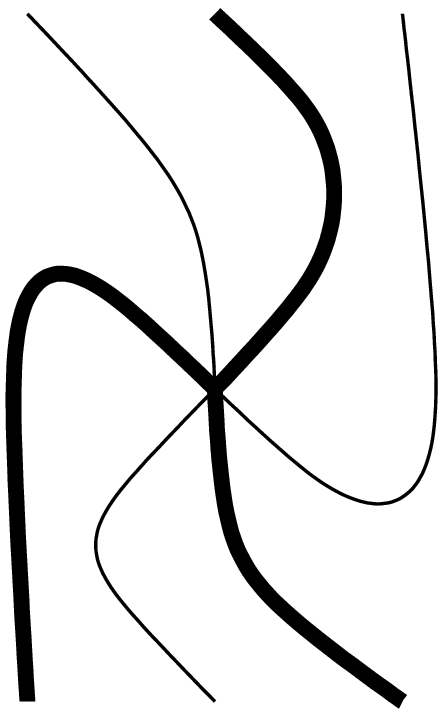} \ \ \  =\ \ \   \ig{.4}{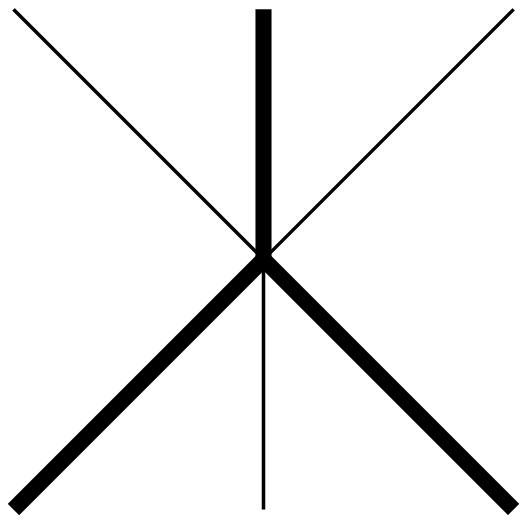} \ \ \ = 
\ \ \   \igh{.32}{ipipipRot1.eps} \label{eqn-ipipipRot}
\end{equation}

Direct sum decompositions of $B_{i,i+1,i}$ and $B_{i+1,i,i+1}$ are encoded by the 
relation 
\begin{equation}
\ig{.3}{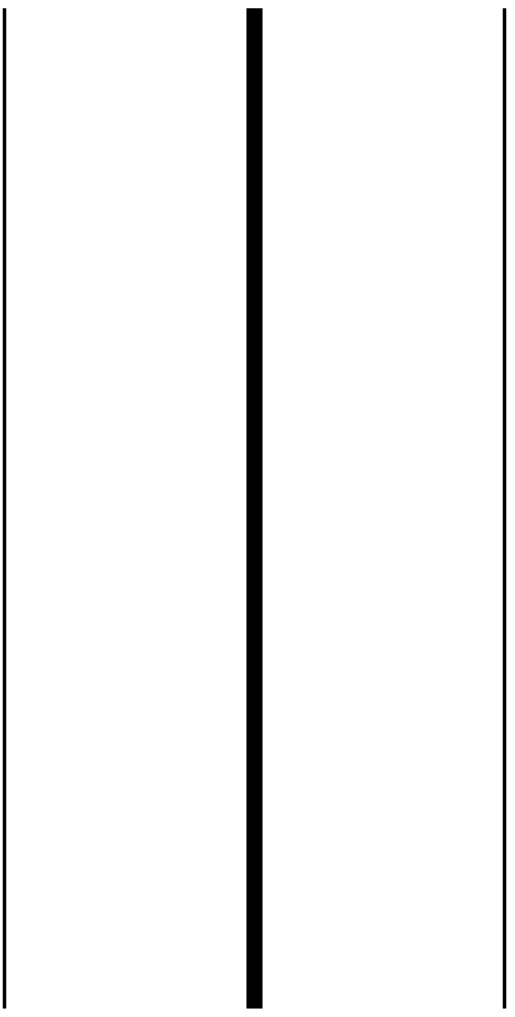}\ \ =\  \ \ig{.3}{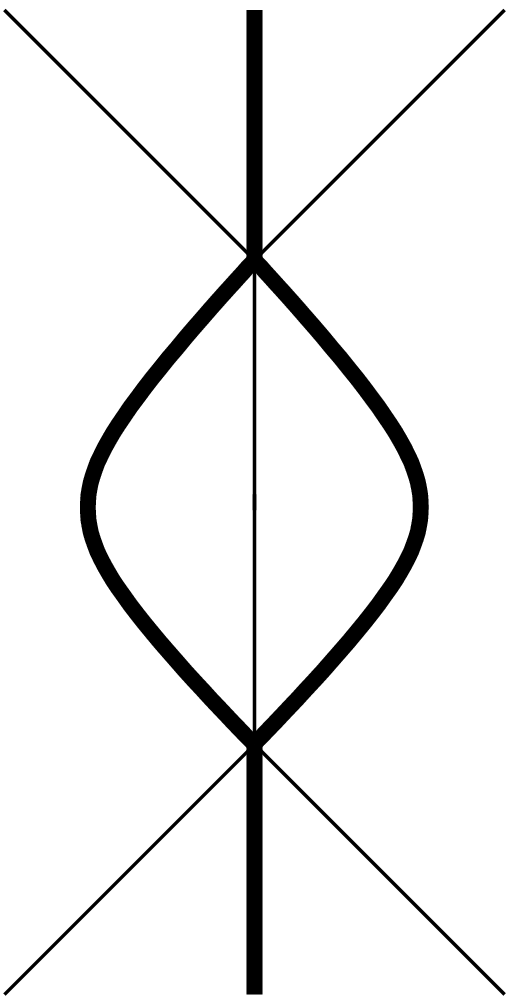} \ \ - \ \ \ig{.3}{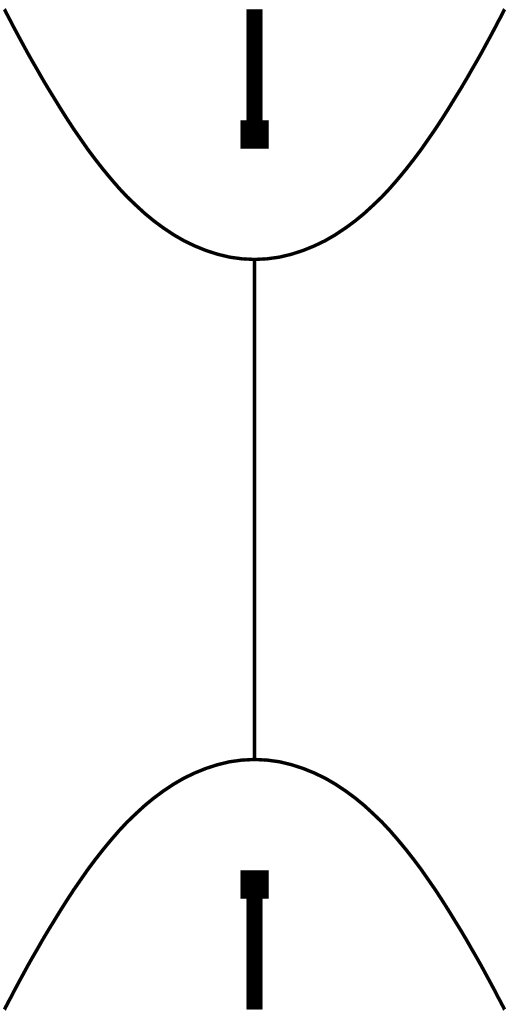} 
\label{eqn-threeLines}
\end{equation} 
and its relative given by reversing the thickness of lines. Other relations are 
(add their reverses as well)  

\begin{equation}
\ig{.3}{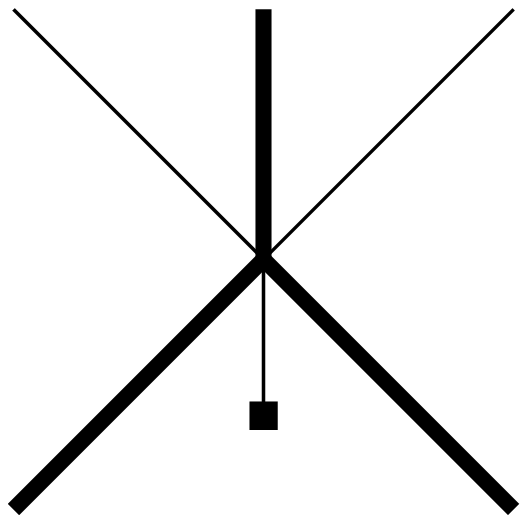} \ \ = \ \ \ \ig{.3}{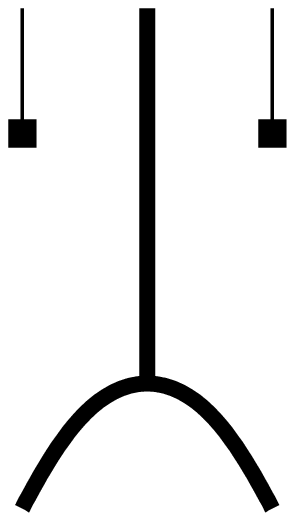} \ \ + \ \ \igh{.3}{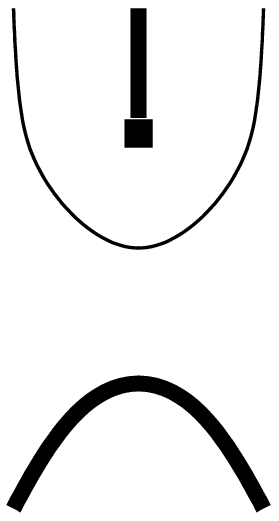} 
\label{eqn-ipipipDot}
\end{equation}

\begin{equation}
\ig{.3}{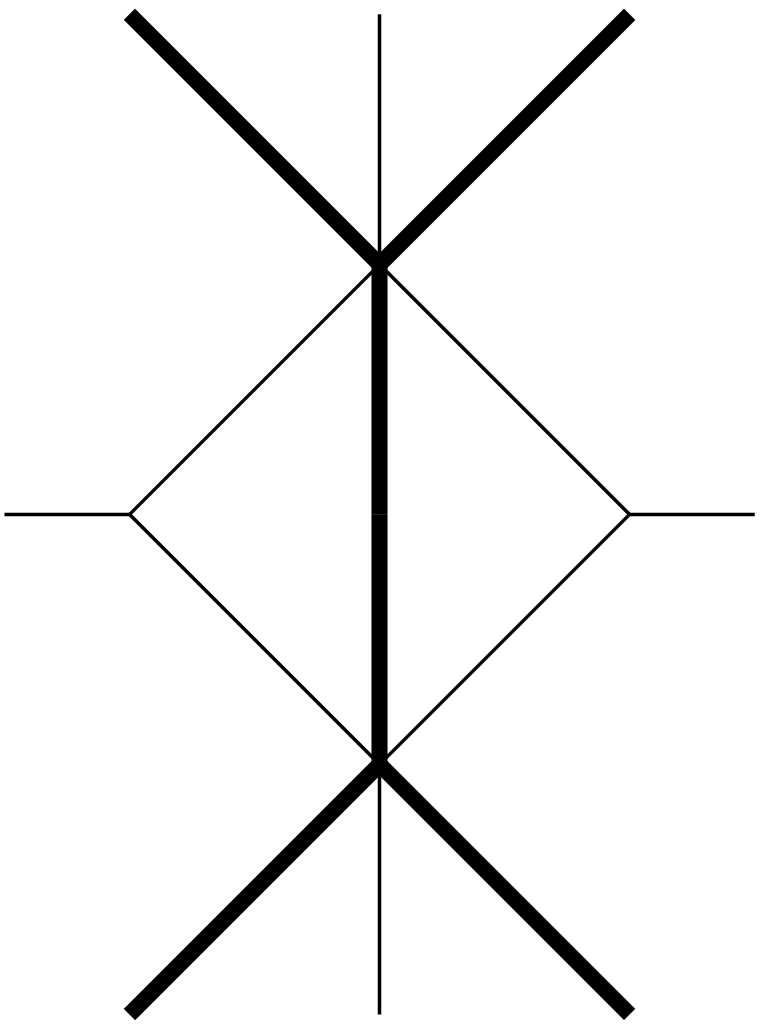} \ \ = \ \ \igrotCW{.3}{ipipipAss.eps} \ \ = 
\ \ \ig{.3}{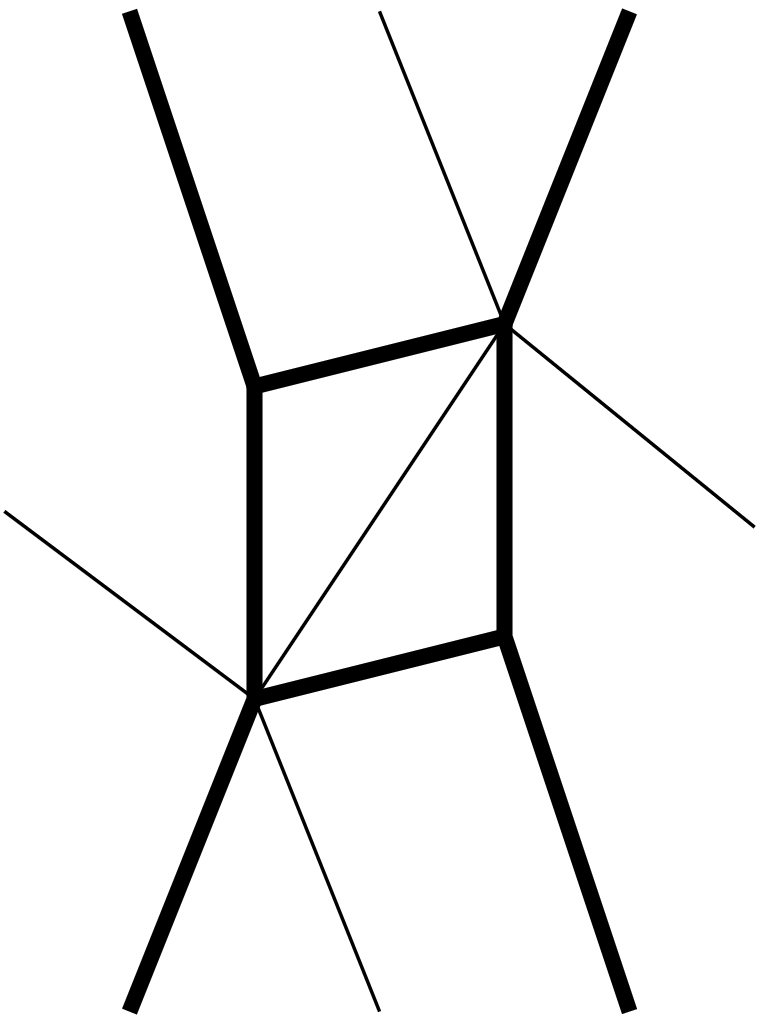} \ \ = \ \  \igh{.3}{ipipipAss2.eps} \label{eqn-ipipipAss}
\end{equation}

\begin{equation}
\ig{.3}{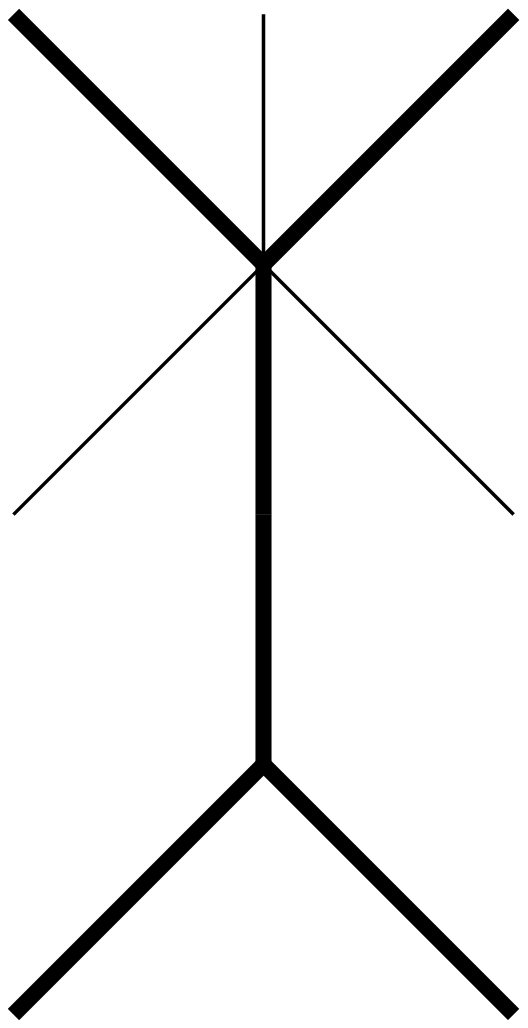} \ \ \ \ = \ \ \ \  \ig{.3}{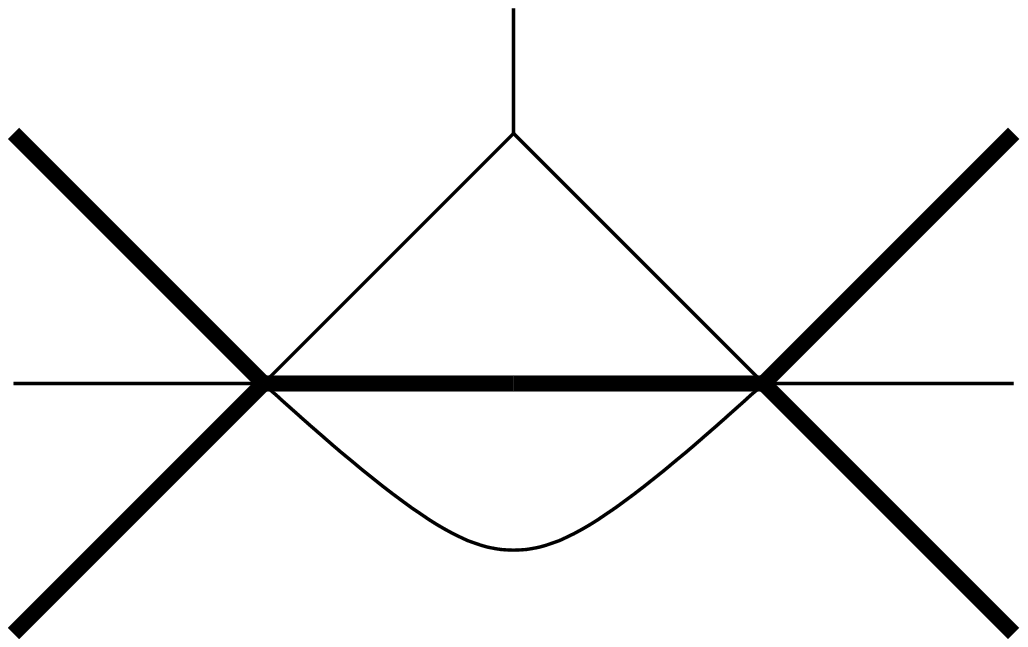} \label{eqn-ipipipAssWDot}
\end{equation}

Lastly, $i$ and $j$-colored lines are allowed to cross when $|i-j|>1$. These crossings 
conclude the list of generating 2-morphisms.  The following relations say that the 
crossings are essentially ``virtual'', i.e. the lines freely go through each other (the solid line 
has color $i$ and dashed -- color $j$). 

\begin{equation}
\ig{.4}{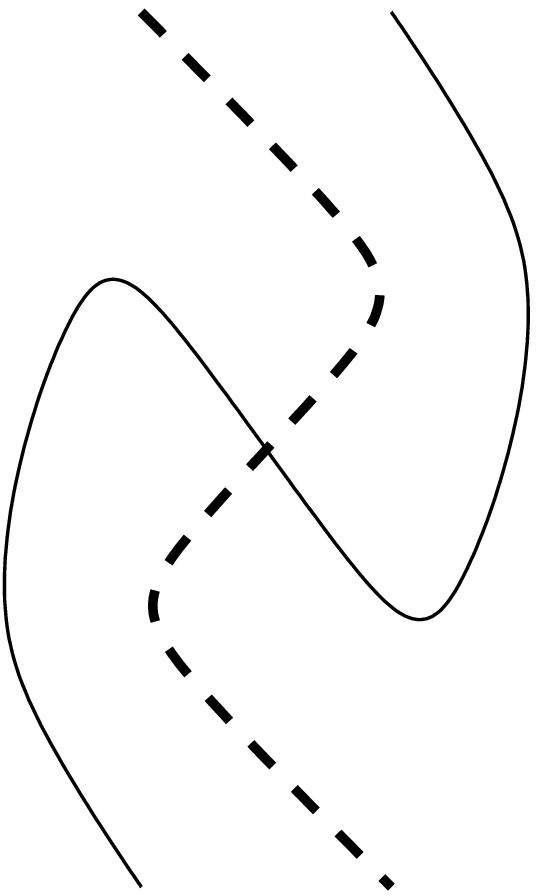}\ \ = \ \ \ig{.4}{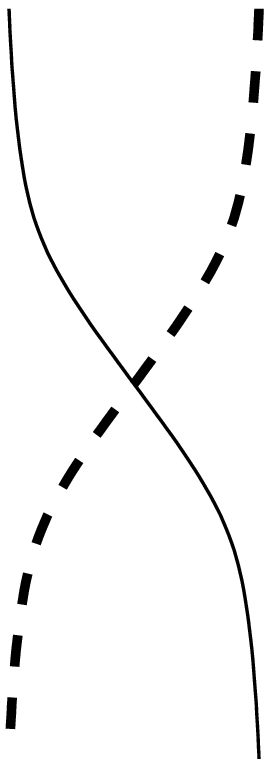}\ \ =\ \ \ig{.4}{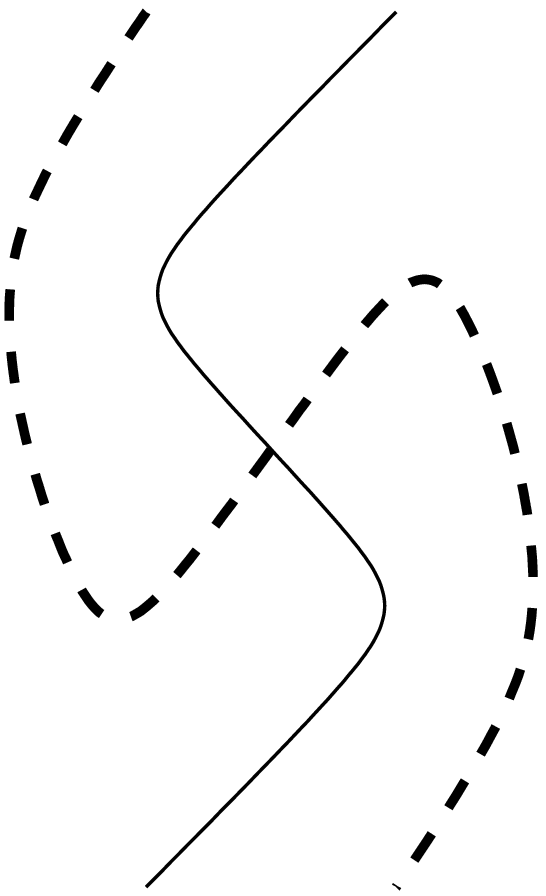} \label{eqn-ijijRot}
\end{equation}
\begin{equation}
\ig{.5}{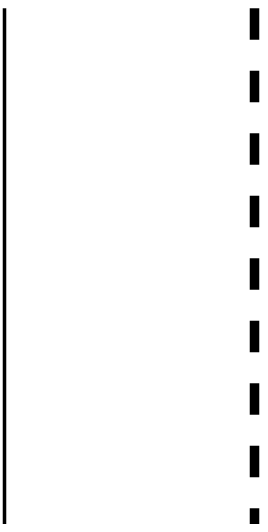}\ \ =\ \ \ig{.5}{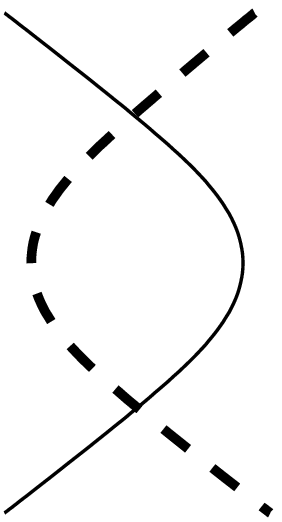} \label{eqn-R2move}  
\end{equation}
\begin{equation}
\ig{.6}{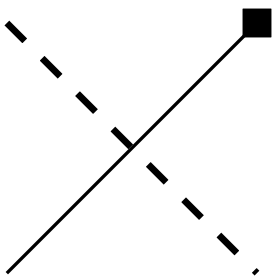} \ \ = \ \ \ig{.6}{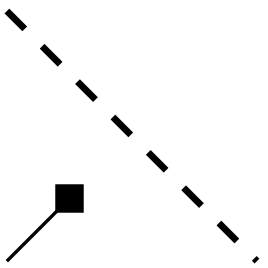} \label{eqn-ijijDot}
\end{equation}
\begin{equation}
\ig{.4}{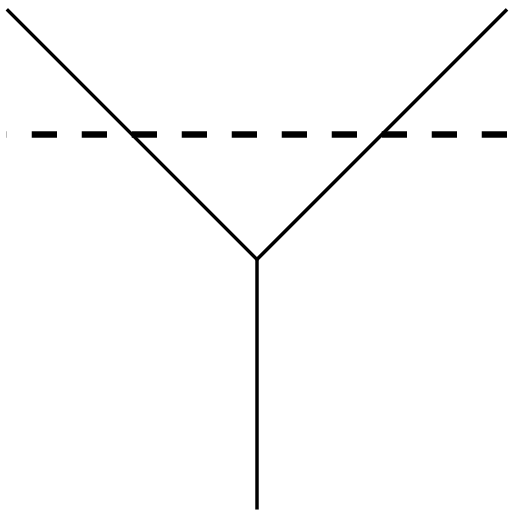} \ \ = \ \ \ig{.4}{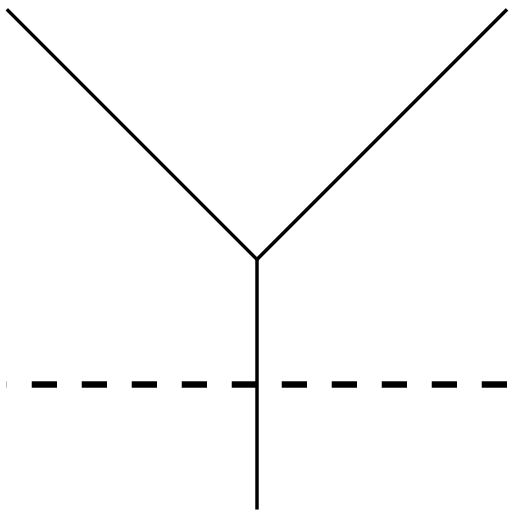} 
\label{eqn-pullFarThruTrivalent} 
\end{equation}

The next relation says that $j$-line interacts trivially with trivalent $i,i+1$-colored 
6-vertex (necessarily $j<i-1$ or $j>i+2$). 
\begin{equation}
\ig{.4}{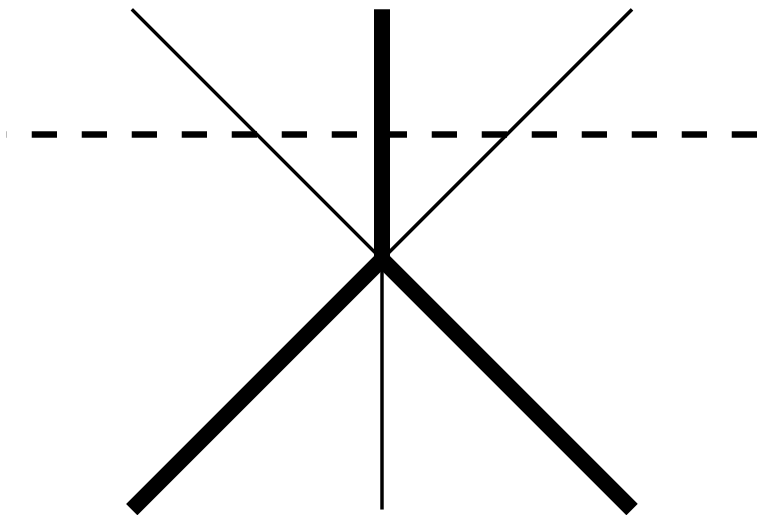} \ \ = \ \ \ig{.4}{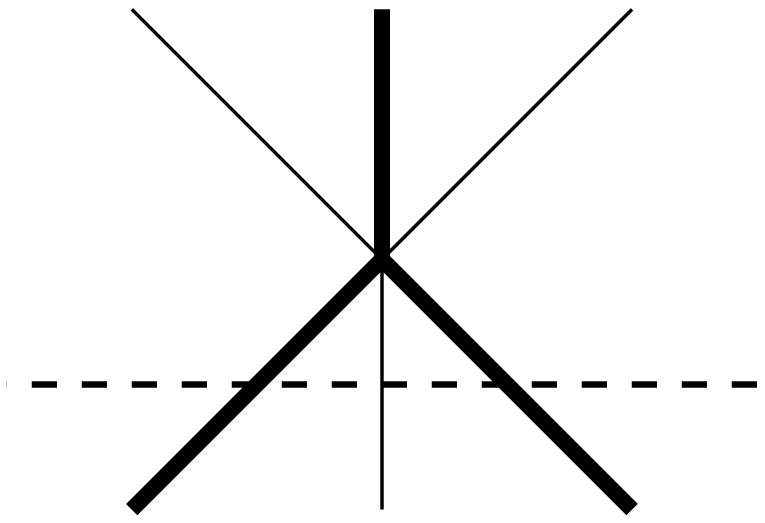}  
\label{eqn-pullFarThru6Valent} \end{equation}
The following relation shows trivial interaction between the crossings of 
$i,j,k$-colored lines (the colors are necessarily far apart $|i-j|>1, |i-k|>1, |j-k|>1$). 
\begin{equation}
\ig{.4}{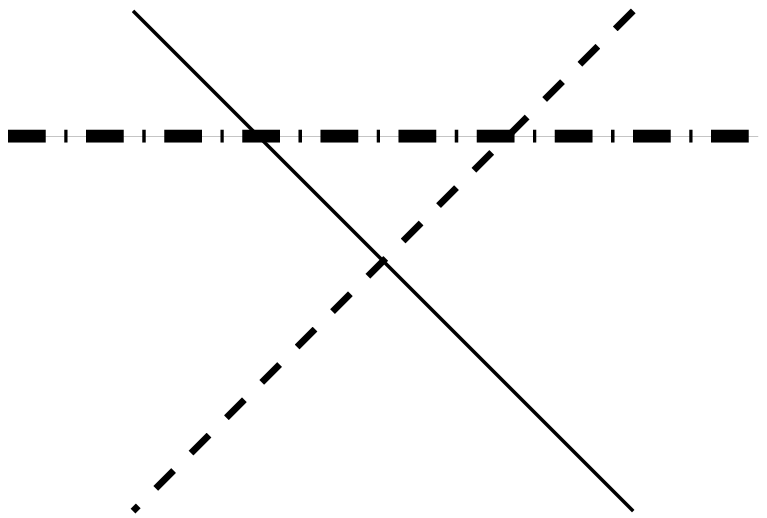} \ \ = \ \ \ig{.4}{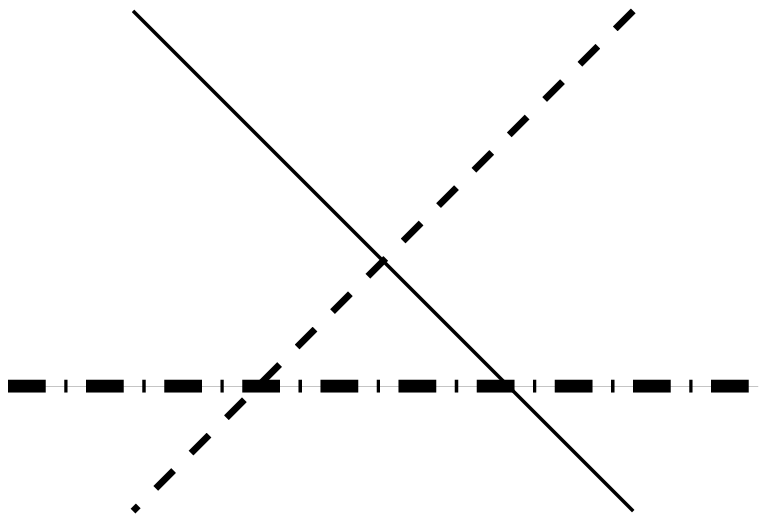}  \label{eqn-R3move}
\end{equation}

Finally, the last and most sophisticated relation is an interation between 
6-valent vertices for $i,i+1,i+2$ (dotted line has color $i+2$).

\begin{equation} \ig{.4}{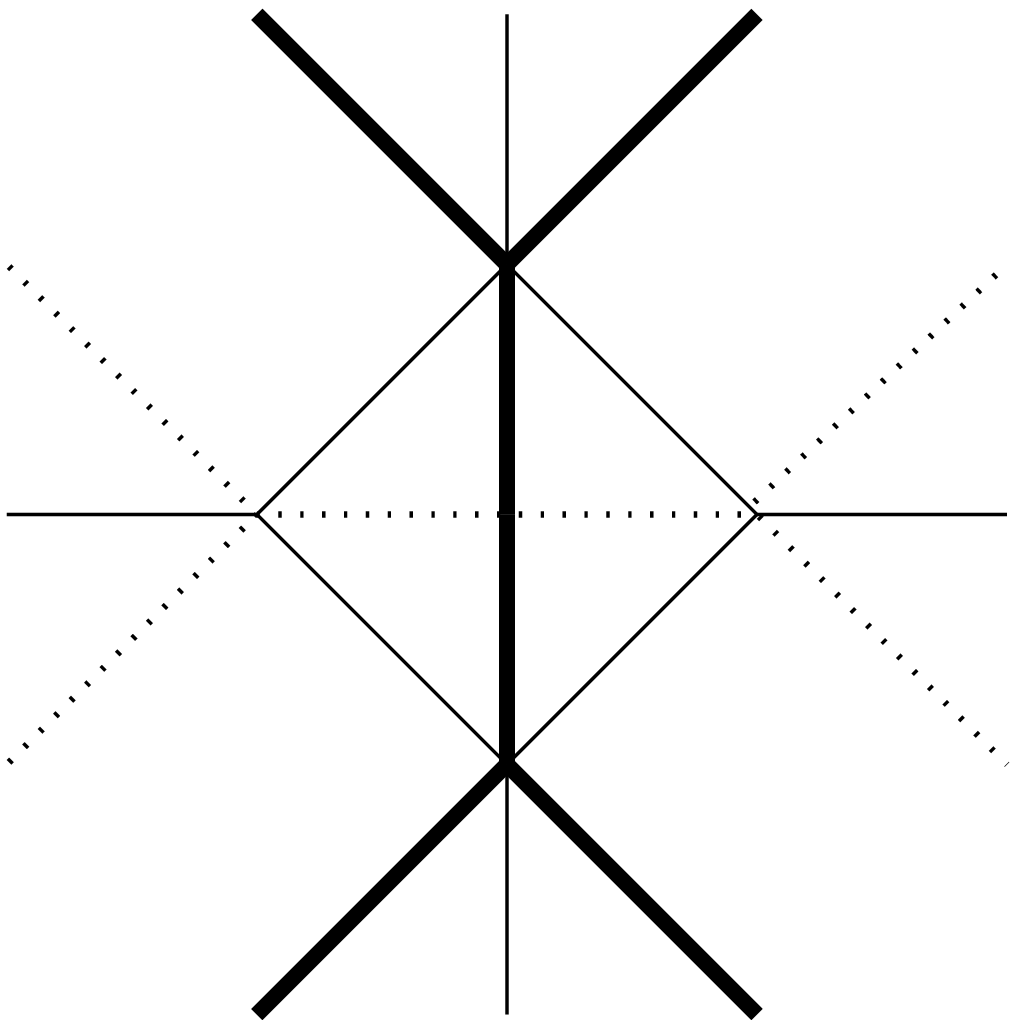} =
  \ig{.4}{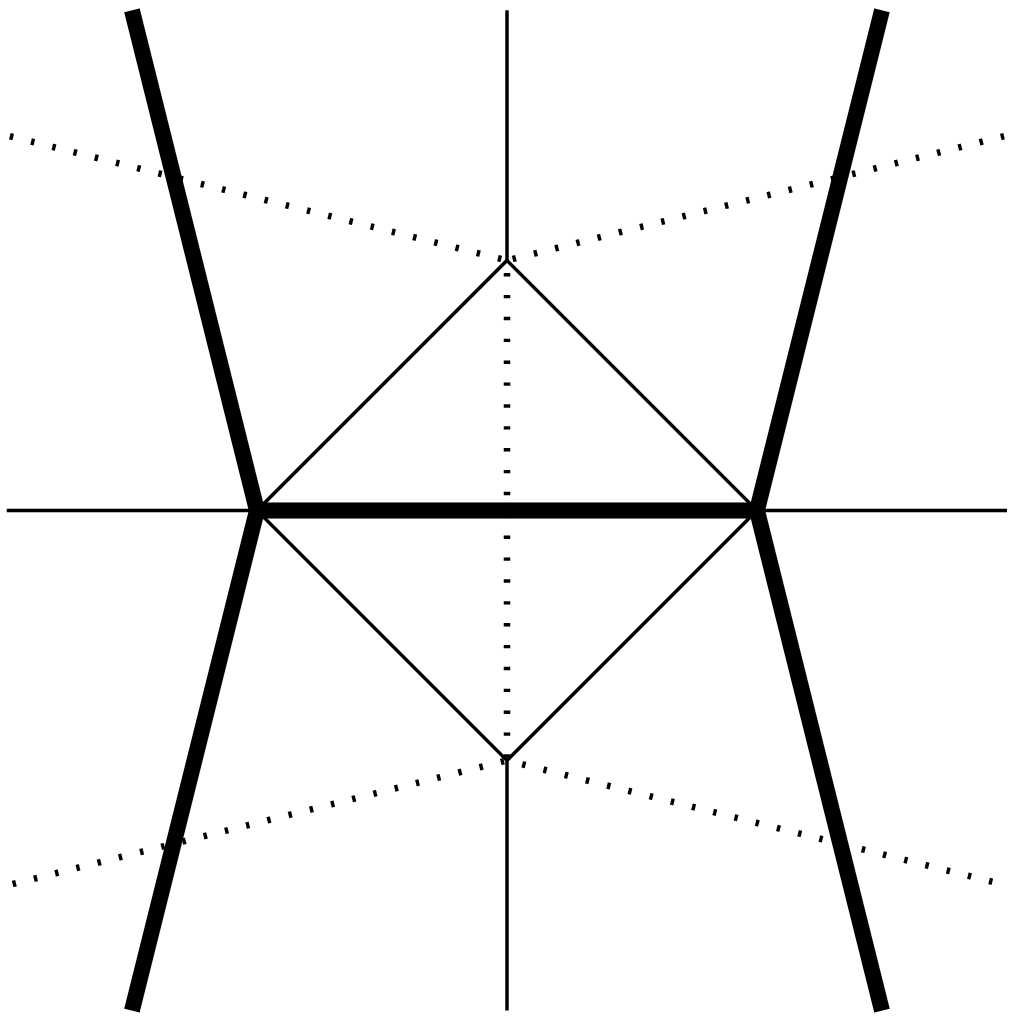} \label{eqn-threeColorAss} \end{equation}

The main result of~\cite{EK} is that the above generators and relations 
give a presentation of the Soergel category $\mc{SC}_1$. We can say 
that $\mc{SC}_1$ is a finitely-presented $\C$-linear pivotal monoidal 
category. This presentation 
is manifestly planar (perhaps the term \emph{planar category} can be used as a sibstitute for 
\emph{pivotal monoidal category}). Of course, the self-adjointness of tensoring with $B_i$ 
was a strong hint that $\mc{SC}_1$ should have a planar description. 

Below we list some applications of this new viewpoint on the Soergel category. 

\vspace{0.1in} 

(a) 
Planar presentation of $\mc{SC}_1$ leads to a new categorification~\cite{E} of the 
Temperley-Lieb algebra $TL_n$.   
The Temperley-Lieb algebra is the quotient of the Hecke algebra $H_n$ by 
the relations $b_i b_{i\pm 1}b_i = b_i$. This algebra is fundamental for the 
construction of the Jones polynomial and admits a categorification via bimodules 
over certain rings, see~\cite{FVIT}.  To categorify the relation $b_i b_{i\pm 1}b_i = b_i$ 
observe that after categorification both sides become Soergel bimodules 
$B_{i,i\pm 1, i}$ and $B_i$, and that the latter bimodule is a summand of the former. 
Thus, the equality can be categorified by setting the complementary summand to $0$. 
Since the 6-valent vertex (\ref{eqn-ipipipRot}) 
is a map which goes through this summand, we simply 
set all 6-valent vertices to 0 and form the quotient monoidal category. A result of 
B.~Elias~\cite{E} says that the Grothendieck ring of the quotient category is 
naturally isomorphic to (a $\Z[q,q^{-1}]$-form) of the Temperley-Lieb algebra. 

In an earlier categorification of the Temperley-Lieb algebra~\cite{FVIT}, a 
basis of homomorphisms 
between bimodules categorifying products of $b_i$'s was given by 3-dimensional 
objects, namely suitable decorated surfaces with boundary and corners embedded 
in $[0,1]^3$. In Elias's approach, a basis is given by a certain collection of 
planar diagrams representing morphisms in the quotient category. A comparison 
between the two categorifications gives an example of \emph{dimensional encoding} 
or \emph{dimensional reduction}, when 3-dimensional information is flattened onto 2D, 
see~\cite{MV, V}. 

\vspace{0.1in} 

(b) R.~Rouquier~\cite{Rou1} constructed a braid group action on the 
category of complexes of Soergel bimodules up to chain homotopies. The braid 
group generators $\sigma_i$ correspond to complexes 
 $$ 0 \lra B_i \{1\} \lra R \lra 0 $$
with the differential being the multiplication map $R\otimes_{R^i} R\lra R$ (counit 
of $B_i$ in graphical notation). This braid group action extends to an action of 
the category of braid cobordisms~\cite{EKr, KT}. The action can be rethought in the 
above diagrammatical language~\cite{EKr}, providing a link between 2-dimensional defining 
relations in $\mc{SC}_1$ and braid cobordisms (4-dimensional objects). 

\vspace{0.1in} 

(c) Taking Hochschild homology of Rouquier complexes in a suitable way 
leads to homology groups that turn out to depend only on the closure 
of a braid~\cite{Ktriply}. The resulting homology theory is triply-graded, 
coincides with the one introduced in~\cite{KR}, and can be 
viewed as a categorification of the Ocneanu trace, thus a categorification of the 
HOMFLY-PT polynomial. We hope that planar diagrammatics will help to 
understand this largely mysterious homology theory and its generalization by 
Webster and  Williamson~\cite{WW}. 

\vspace{0.1in} 

It is natural to expect a possible relation between 
diagrammatics for $\mc{SC}_1$ (and more general planar categories and 
planar 2-categories) and 2D (topological) field theories. 
Biadjointness is implicit everywhere in 2D QFTs, 
since the categories that appear there, such as derived 
categories of coherent sheaves on 3D Calabi-Yau manifolds, Fukaya-Floer 
categories, and categories of matrix factorizations in LG models admit a 
wealth of biadjoint functors, that come up in a natural way, via convolutions 
with sheaves on products of a pair of Calabi-Yau's, via convolutions with 
Lagragians in the product of sympectic manifolds, and via matrix factorizations 
with potentials $f(x)-g(y)$. Planar diagrams for morphisms in $\mc{SC}_1$ (and 
planar diagrams for 2-morphisms in similar 2-categories)  
can perhaps be viewed as world sheets of 2D QFTs with defect lines and vertices, 
with regions of the diagram labelled by different target objects.  

On a speculative note, we suggest a noncommutative geometry  
language for $\mc{SC}$. Recall that we represent Soergel bimodule 
$B_{\ii}$ for a sequence $\ii=i_1\dots i_d$ 
 by placing dots labelled $i_1, \dots, i_d$ on a line. Let's imagine that this 
line with labeled dots is a path, and refer to $B_{\ii}$ as a path as well. More generally, 
we call any object of $\mc{SC}$ a path. Indecomposable objects $B_w$  
are \emph{geodesics}. Given two objects $M,N$ of $\mc{SC}$, we think of 
the graded vector space $\Hom_{\mc{SC}}(M,N)$ as 
categorified quantized area of surfaces stretched between $M$ and $N$. 
 The fact that the grading of  $\Hom_{\mc{SC}}(M,N)$ is bounded from 
below is loosely analogous to the same property of energy.


\section{Categorification of quantum groups}

The first example of a categorification of a bialgebra is apparently due to 
L.~Geissinger~\cite{Gei}. 
Standard inclusions of symmetric groups $S_n\times S_m \subset S_{n+m}$ give 
rise to induction and restriction functors 
\begin{eqnarray*} 
\mathrm{Ind}_{n,m}  & : & \C[S_n]\dmod \times \C[S_m]\dmod \lra \C[S_{n+m}]\dmod, \\
  \mathrm{Res}_{n,m} & : & \C[S_{n+m}]\dmod \lra \C[S_n \times S_m]\dmod 
\end{eqnarray*}
 between categories of (finite-dimensional) modules over these group algebras. 
Summing over all $n,m\ge 0$ produces functors 
\begin{equation*} 
\mathrm{Ind} \ : \  \mc{S} \times \mc{S} \lra \mc{S} \ \ \ \ \ \ \ \ \ 
  \mathrm{Res} \ : \ \mc{S} \lra \mc{S} \otimes \mc{S}, 
\end{equation*} 
where 
$$\mc{S} : = \oplusop{n\ge 0} \C[S_n]\dmod, \ \ \ \ \ \ \ \ \ \ 
   \mc{S} \otimes \mc{S} := \oplusop{n,m\ge 0} \C[S_n\times S_m]\dmod
$$ 
is the direct sum of representation categories of symmetric groups, respectively 
products of symmetric groups. 
Functors $\mathrm{Ind}$ and $ \mathrm{Res}$ induce homomorphisms of 
Grothendieck groups 
\begin{equation*} 
[\mathrm{Ind}]  \ : \ K_0(\mc{S}) \otimes K_0(\mc{S}) \lra K_0(\mc{S}), \ \ \ \ \ \ 
  [\mathrm{Res} ]  \ : \ K_0(\mc{S}) \lra K_0(\mc{S}) \otimes K_0(\mc{S}).  
\end{equation*} 
Here 
$$ K_0(\mc{S}) := \oplusop{n\ge 0} K_0(\C[S_n]\dmod),$$
and the Grothendieck group 
$K_0(\C[S_n]\dmod)$ is a free abelian group with a natural basis given 
by symbols of simple $S_n$ modules $L_{\lambda}$, parametrized by 
partitions $\lambda$ of $n$ (the group algebra $\C[G]$ is semisimple 
for any finite group $G$, any module is projective, and isomorphism classes 
of simple modules give rise to a basis of $K_0$). It is natural to identify 
$K_0(\mc{S})$ with the ring of symmetric functions in countably many variables by 
taking $[L_{\lambda}]$ to the Schur function $s_{\lambda}$. 

\begin{prop} Maps $[\mathrm{Ind}]$ and $[\mathrm{Res}]$ turn $K_0(\mc{S})$ 
into a biring. 
\end{prop} 

A biring is a bialgebra over $\Z$. 
This result, due to Geissinger~\cite{Gei}, also tells us that the ring of symmetric functions 
has a natural comultiplication. 
The condition that comultiplication $\Delta=[\mathrm{Res}]$ is a homomorphism, 
$\Delta(xy)= \Delta(x)\Delta(y)$, follows from the Mackey's decomposition theorem 
specialized to induction and restriction between products of symmetric groups. 

Deep generalizations of this construction were developed by A.~Zelevinsky~\cite{Zel}, who 
constructed a similar biring structure on the sum (over all $n$) of Grothendieck groups of 
representation categories of $GL(n,\mathbb{F}_q)$, with parabolic induction and 
restriction in place of the usual induction and restriction. Many other examples of a biring 
structure on Grothendieck groups can be found in Bergeron and Li~\cite{BL} and 
references therein. 

In another paper, Zelevinsky classified irreducible representations of the affine Hecke algebra
 of $S_n$ in terms of multisegments~\cite{Zel2}. 
The significance of his result became clear much later, 
when S.~Ariki~\cite{Ar} categorified all integrable irreducible representations 
of $\mf{sl}_k$ and 
affine $\mf{sl}_k$ via blocks of representation categories of cyclotomic quotients of 
affine Hecke algebra (an important earlier milestone was the work of Lascoux-Leclerc-Thibon 
on categorification of level one representations of affine $\mf{sl}_k$ via the representation 
categories of Hecke algebras of $S_n$ at $k$-th root of unity, sum over all $n\ge 0$).  
I.~Grojnowski~\cite{Gr}, armed with ideas of A.~Kleshchev from modular 
representation theory, gave an alternative  derivation of Ariki's results. 

Ariki also showed that Grothendieck groups of completions of affine Hecke algebra at 
suitable central ideals are canonically isomorphic to weight spaces of $U^+(\mf{sl}_k)$ and 
$U^+(\hat{\mf{sl}}_k)$, giving a conceptual categorification-style explanation for 
Zelevinsky's classification of irreducibles and proved that the basis of indecomposable 
projectives in these Grothendieck groups coincides with the $q=1$ 
specialization of the canonical basis in the quantum group $U^+_q$. 

In Lusztig's construction of the canonical basis~\cite{Lu, Lus2}, 
partially inspired by Ringel's work~\cite{Rin}, canonical basis elements correspond to 
simple perverse sheaves on the varieties of quiver representations.  
In fact, his construction can even be viewed as a categorification of $U^+_q$ if one passes 
from the set of simple perverse sheaves to a suitable triangulated category of 
equivariant sheaves that these simple objects generate. Lusztig~\cite{Lu} and 
Kashiwara~\cite{Ka} also provided more elementary approaches to the 
canonical basis. 

\vspace{0.06in}

We'll now review a very down-to-earth categorification of $U^+_q$ based on a collection 
of certain rings, following~\cite{KL1}, see also \cite{Rou2}. 
Grothendieck groups of these rings can be identified with an integral 
version of $U^+_q$, while induction and restriction functors between these rings 
correspond 
to multiplication and comultiplication in the quantum group. A direct link 
between this categorification of $U^+_q$ and Lusztig's perverse sheaves construction was 
found by Varagnolo and Vasserot~\cite{VV}. Brundan and Kleshchev~\cite{BK} 
related these rings to Ariki's categorifications of $U^+$ and highest weight 
representations. A direct relation to derived categories of
coherent sheaves on Nakajima quiver varieties~\cite{Nak} is expected~\cite{CK}. 

\vspace{0.2in} 

Let $\Gamma$ be an unoriented graph without loop and multiple edges, with set of 
vertices $I$. The quantum group $U^+=U^+(\Gamma)$ 
is a $\Z[I]$-graded $\Q(q)$-algebra with generators $E_i, i\in I$ of degree $i$ and 
defining relations 
\begin{eqnarray} 
   E_i E_j & = & E_j E_i \ \  \mbox{if $(i,j)$ is not an edge}, \label{deq1}\\
   (q+q^{-1})E_i E_j E_i  & = & E_i^2 E_j + E_j E_i^2 \ \  \mbox{if $(i,j)$ is an edge}. 
 \label{deq2} 
\end{eqnarray} 
The integral version $U^+_{\Z}$ is the $\Z[q,q^{-1}]$-subalgebra of $U^+$ 
generated by divided powers $E_i^{(n)}:= \frac{E_i^n}{[n]!}$ over all $i$ and $n$. 

Both $U^+$ and $U^+_{\Z}$ are twisted bialgebras (the Cartan 
subalgebra is missing), with $\Delta(E_i)=E_i\otimes 1 + 1\otimes E_i$
(instead of the usual comultiplication $\Delta(E_i)=E_i\otimes 1 + K_i\otimes E_i$
in quantum groups), and 
nonstandard algebra structure on $U^+\otimes U^+$: 
 \begin{eqnarray}
(x_1\otimes x_2) (x_1'\otimes x_2') = q^{ -|x_2|\cdot |x_1'|} x_1 x_1' \otimes
x_2 x_2'
\end{eqnarray}
where $|x|$ is the degree, taking values in $\Z[I]$, and $\cdot$ is the inner product 
on $\Z[I]$ with $i\cdot i =2$, $i\cdot j=-1$ if $i$ and $j$ are connected by an edge, 
and $i\cdot j=0$ otherwise. Relative to this algebra structure, $\Delta$ is a 
homomorphism $U^+\lra U^+\otimes U^+$ (see~\cite{Lus2} for more details).  

Let $\nu\in \N[I], \ \  \nu = \sum_{i \in I} \nu_{i} \cdot i \;, \ \  \nu_i \in \N.$ 
For each such $\nu$ define the graded ring $R(\nu)$, as ring spanned  by 
diagrams of lines in $\R\times [0,1]$, with $\nu_i$ lines colored by label $i\in I$. 
Lines can intersect, but triple intersections are not allowed, neither are U-turns 
that create critical points under projection onto the $y$-axis. Lines may carry dots
 (in the picture below $\nu=2i+j+k$). 

\begin{center}{\psfig{figure=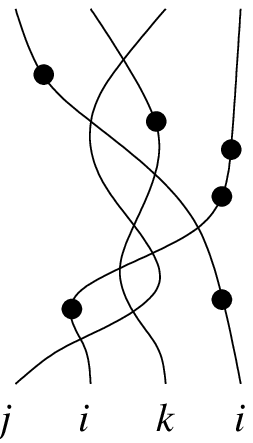,height=4cm}}\end{center}

Product is given by concatenation of diagrams; if the labels of endpoints don't match 
the product is zero. 
We allow isotopies (that do not create critical points relative to the $y$-axis projection) 
and impose the following relations.  

\begin{eqnarray*} 
& & {\psfig{figure=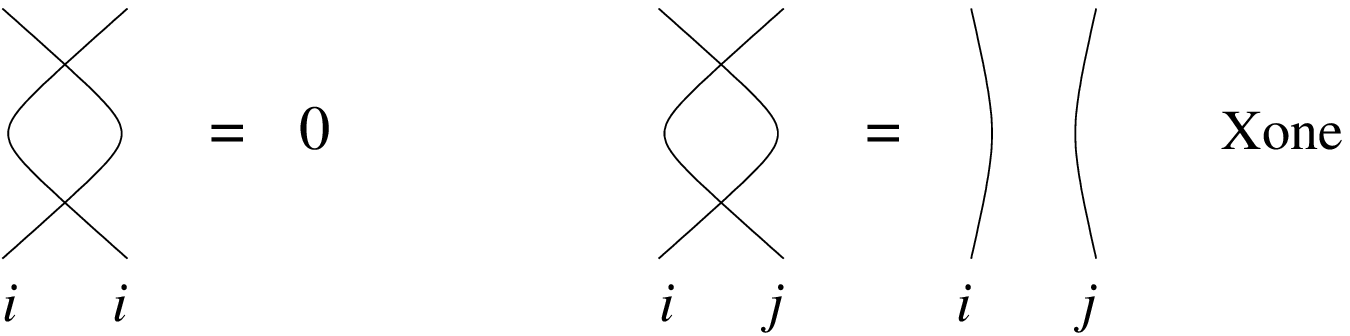,height=2.5cm}}    \\
& & \quad  \\ 
& & {\psfig{figure=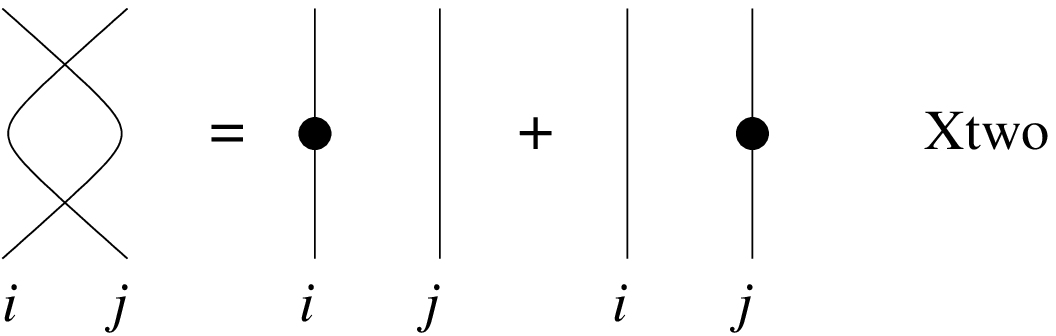,height=2.5cm}} 
\end{eqnarray*} 

\begin{equation*} 
{\psfig{figure=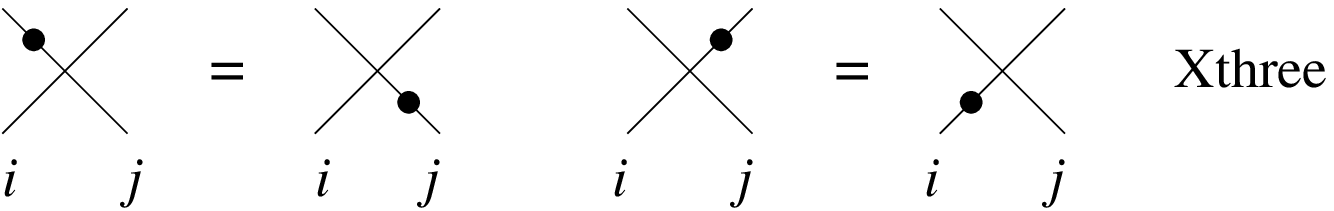,height=2cm}}  
\end{equation*} 

\begin{equation*} 
{\psfig{figure=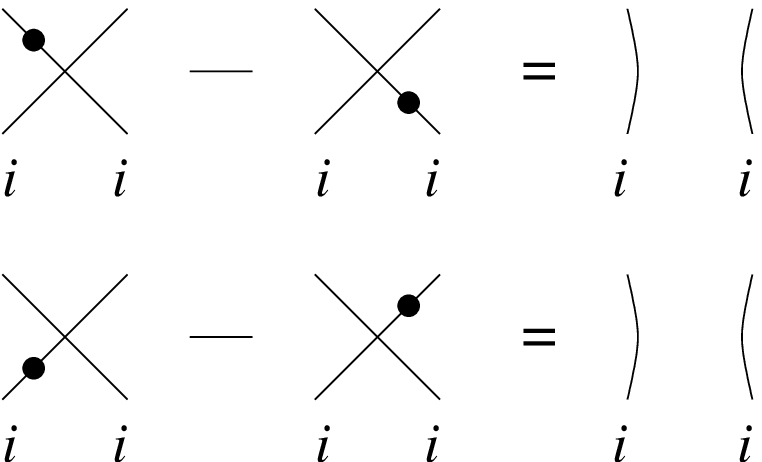,height=4cm}}   
\end{equation*} 

\begin{equation*} 
{\psfig{figure=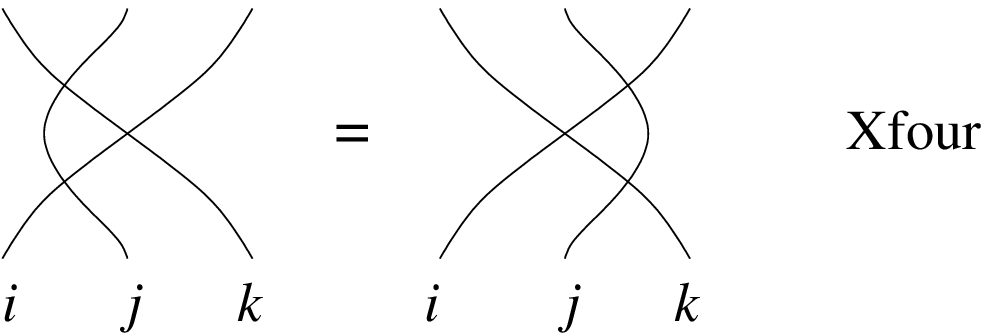,height=3cm}}   
\end{equation*}

\begin{equation*} 
{\psfig{figure=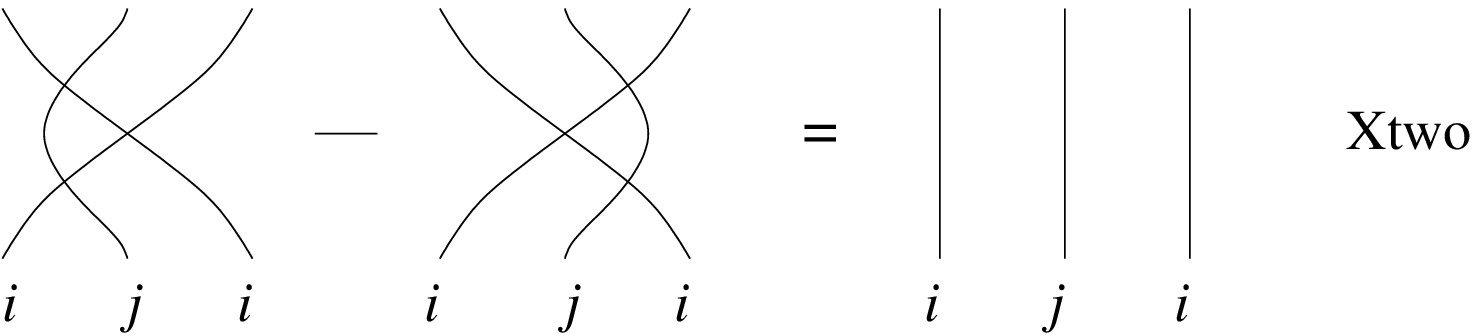,height=3cm}}   
\end{equation*}

Make $R(\nu)$ graded by assigning degree $2$ to a dot and degree $-i\cdot j$ to 
the $(i,j)$-crossing. Also, choose a base field $\Bbbk$ and consider $R(\nu)$ as 
a graded $\Bbbk$-algebra. Diagrammatics make it clear that, at the cost of 
adding simpler diagrams, all dots can be moved above all crossings. Moreover, 
if two lines intersect twice in a diagram, the diagram can be further simplified. 
This argument gives a spanning set in $R(\nu)$ which consists of 
all monomials in dots times minimal length presentations of symmetric group 
elements as products of crossings, with arbitrary $\nu$-colorings of the lines. 
One can check that this spanning set is a basis, by looking at the action of $R(\nu)$ 
on a suitable representation. 

Each sequence $\ii = i_1\dots i_m$ of weight $\nu$ gives rise to an idempotent 
$1_{\ii}$ in $R_{\nu}$ given by the diagram with $m$ vertical lines labelled 
$i_1, \dots, i_m$ from left to right. These idempotents are mutually orthogonal 
and $1=\sum_{\ii\in Seq(\nu)} 1_{\ii}$, where $Seq(\nu)$ is the set of sequences 
of weight $\nu$. Also, each idempotent $1_{\ii}$ determines a graded 
projective module $P_{\ii} =R(\nu) 1_{\ii}$. 

\begin{prop} \cite{KL1, Rou2} There are natural isomorphisms of graded projectives 
\begin{eqnarray*} 
 P_{ij}  & \cong & P_{ji} \hspace{0.2in} \mbox{if} \hspace{0.1in} i\cdot j=0, \\
 P_{iji}\{1\} \oplus P_{iji}\{-1\} & \cong & P_{iij} \oplus P_{jii} 
 \hspace{0.1in} \mbox{if} \hspace{0.1in} i\cdot j=-1. 
\end{eqnarray*} 
\end{prop} 

This is a crucial proposition. After passage to the Grothendieck 
group, these isomorphisms become equalities on symbols of projectives: 
\begin{eqnarray*} 
 [P_{ij}]  & =  & [P_{ji}] \hspace{0.2in} \mbox{if} \hspace{0.1in} i\cdot j=0, \\
 q [P_{iji}] +  q^{-1} [P_{iji}] & =  & [P_{iij}] +  [P_{jii}] 
 \hspace{0.1in} \mbox{if} \hspace{0.1in} i\cdot j=-1,  
\end{eqnarray*} 
showing a perfect match with the defining relations (\ref{deq1}), (\ref{deq2}) in $U^+$
after converting $[P_{ij}]$ to $E_i E_j$, $[P_{iji}]$ to $E_i E_j E_i$, etc. 

The second isomorphism can be refined. Specifically, 
$R(m i)$ ($m$ strands identically labelled) is naturally isomorphic to the 
nilHecke algebra, generated by multiplication by monomials $x_1, \dots, x_m$ 
and divided difference operators $\partial_i$. It is well-known that the nilHecke 
algebra is isomorphic to the algebra of $m!\times m!$-matrices with coefficients 
in the ring of symmetric functions in $x_1, \dots, x_m$. Consequently, 
the free module $R( mi)$ decomposes as a sum of $m!$ copies of an 
indecomposable projective module (the column module) 
which we denote by $P_{i^{(m)}}$. The latter has a suitably selected overall 
grading shift, so that, taking grading into  account, 
 $R( mi)$ is the sum of $[m]!$ copies of $P_{i^{(m)}}$. 

When $m=2,$ there is a decomposition $P_{ii} \cong  P_{i^{(2)}}\{1\} 
\oplus P_{i^{(2)}}\{-1\}$, and the second isomorphism above simplifies to 
 $$ P_{iji}\cong P_{i^{(2)}j} \oplus P_{j i^{(2)}} .$$ 

Just like with the group algebras of symmetric groups, we can fit 
rings $R(\nu)$ over all $\nu\in \N[I]$ into a tower of algebras. 
There are natural (nonunital) inclusions $R(\nu)\otimes R(\nu') \subset R(\nu+\nu')$ 
described by placing diagrams representing elements of $R(\nu)$ and 
$R(\nu')$ in parallel, next to each other. Consider associated induction and 
restriction functors $\Ind_{\nu,\nu'}$, $\Res_{\nu,\nu'}$. We would like to 
look at the maps these functors induce on the Grothendieck group. 
The induction functor always takes projectives to projectives and induces 
a homomorphism between Grothendieck groups of finitely-generated graded 
projective modules 
$$ [\Ind_{\nu,\nu'}]: K_0(R(\nu))\otimes_{\Z[q,q^{-1}]} K_0(R(\nu')) \lra 
 K_0(R(\nu)\otimes_{\Bbbk} R(\nu'))  \lra K_0(R(\nu+\nu')).$$ 
The first arrow comes from tensoring projectives. In our case, the first arrow is 
an isomorphism, due to absolute irreducibility of simple graded 
$R(\nu)$-modules, which we prove using methods of Kleshchev and Grojnowski. 
The second arrow comes from the induction functor. 

For a general inclusion of rings, restriction functor does not necessarily take 
projectives to projectives. However, $R(\nu+\nu')$ (more accurately, 
$(1_{\nu}\otimes 1_{\nu'})R(\nu+\nu')$) is a projective 
$R(\nu)\otimes R(\nu')$-module, so that the restriction does take projectives 
to projectives and induces a homomorphism 
$$ [\Res_{\nu,\nu'}]:  K_0(R(\nu+\nu')) \lra K_0(R(\nu)\otimes_{\Bbbk} R(\nu')) 
\cong K_0(R(\nu))\otimes_{\Z[q,q^{-1}]} K_0(R(\nu')).$$ 
Form the direct sum 
$$ K_0(R) := \oplusop{\nu\in \N[I]} K_0(R(\nu))$$ 
and the corresponding sums of functors 
$$ \Ind:= \oplusop{\nu,\nu'} \Ind_{\nu,\nu'}, \ \ \    
\Res:= \oplusop{\nu,\nu'} \Res_{\nu,\nu'}. $$
These induce homomorphisms of $K_0$-groups 
$$ [\Ind] \ : \ K_0(R)\otimes K_0(R) \lra K_0(R), \ \  \ 
   [\Res] \ : \ K_0(R) \lra K_0(R)\otimes K_0(R).$$ 

\begin{theorem} \cite{KL1} 
There is a canonical isomorphism of twisted bialgebras over $\Z[q,q^{-1}]$ 
 $$ U^+_{\Z}\cong K_0(R)$$ 
under which multiplication, respectively comultiplication, in $U^+_{\Z}$ is 
given by $[\Ind]$, respectively $[\Res]$. 
\end{theorem} 

This isomorphism takes the product element $E_{\ii} = E_{i_1}\dots E_{i_m}$ to 
the symbol $[P_{\ii}]$ of $P_{\ii} = P_{i_1\dots i_m}$ 
and divided power element $E_{i}^{(m)}$ to the symbol 
of $P_{i^{(m)}}$, which is an indecomposable direct summand of $P_{i^m}$ with 
a suitably normalized grading. 

It is shown by Brundan  and Kleshchev~\cite{BK2} 
(in the $sl(k)$ and affine $sl(k)$ case) 
and Varagnolo and Vasserot~\cite{VV} (for general graphs $\Gamma$) that, if 
$\mbox{char}(\Bbbk)=0$,  
the canonical basis~\cite{Lu, Ka, Lus2} of $U^+_{\Z}$ goes to the basis of 
indecomposable projectives under this isomorphism. To get this match, when $\Gamma$ 
has an odd length cycle, the definition of $R(\nu)$ should be slightly modified, 
see the above references and~\cite{KL2, Rou2}. 
Rings $R(\nu)$ proved handy in recent constructions~\cite{BK, BKW, Ar2} 
of graded versions of the group algebra of $S_n$, Specht modules, and $q$-Schur 
algebras, also see~\cite{HMM} and references therein for related developments.  
Moreover, they serve as a building block in Webster's categorification of 
Reshetikhin-Turaev link and tangle invariants~\cite{W1, W2}. 

\vspace{0.1in} 

$U^+$, its categorification, and related structures can be incorporated 
into the following (rather incomplete) diagram, where vertical up arrows denote 
categorifications. 

\begin{center}{\psfig{figure=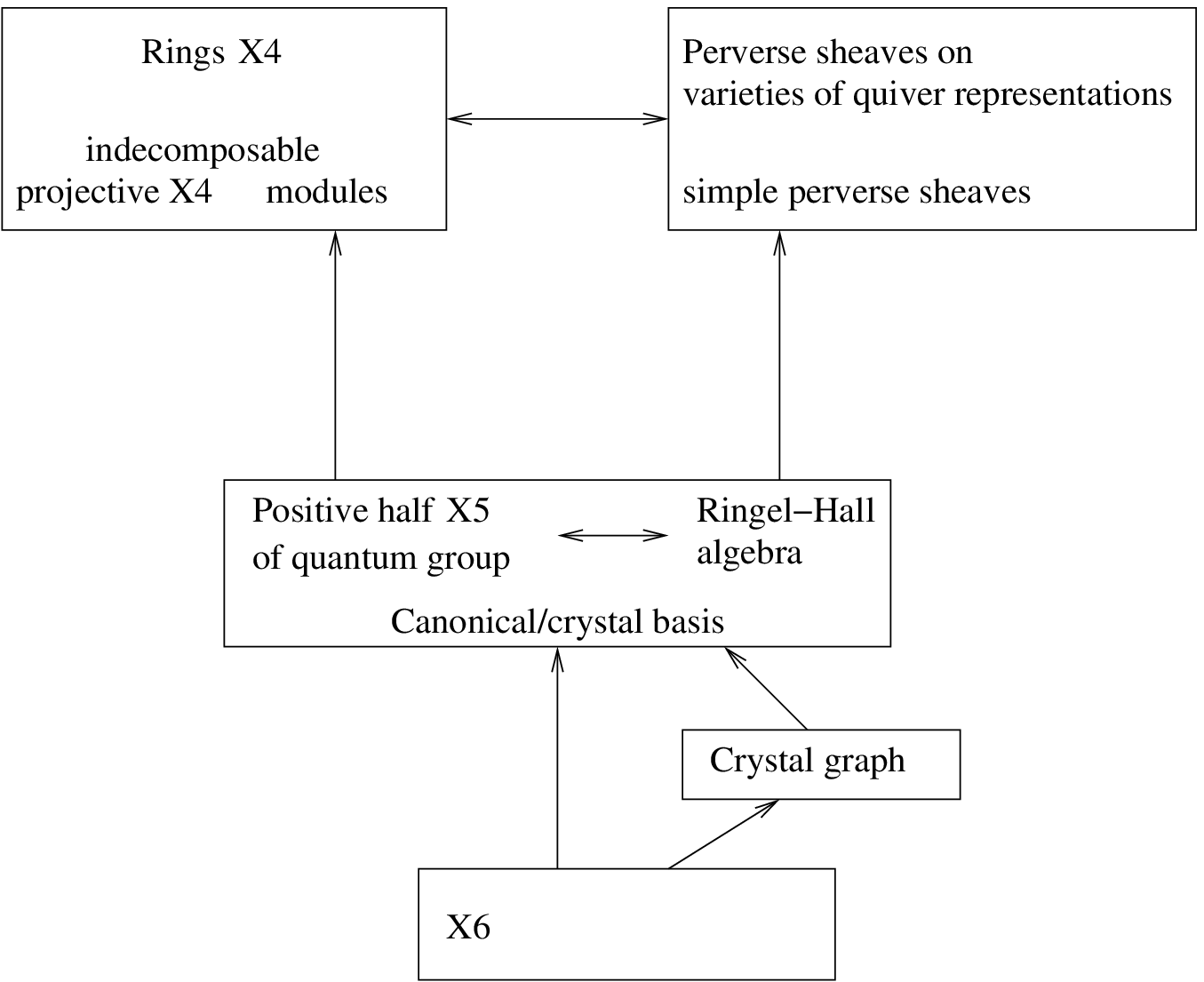,height=10cm}}\end{center}

At the base lies the collection of positive integers $\dim(U^+(\nu))$ -- these 
are dimensions 
of weight spaces $U^+(\nu)$ of $U^+$, or, more implicitly, coefficients 
of the Kostant partition function. Its categorification is the twisted bialgebra 
$U^+$ equipped with the Lusztig canonical/Kashiwara crystal basis. 
After the second round of categorification, the weight spaces $U^+(\nu)$ 
become Grothendieck groups of graded rings $R(\nu)$, and the canonical 
basis elements lift to indecomposable projective modules. Equivalently, 
one can work with the derived category of equivariant constructible sheaves 
on varieties of quiver representations, with simple perverse sheaves 
being analogues of indecomposable projective graded $R(\nu)$-modules. 
The upper horizontal arrow is a derived equivalence (for a carefully chosen 
version of the sheaves category) exchanging 
indecomposable projective $R(\nu)$-modules and simple perverse sheaves. 

Kashiwara crystal graph, together with Kashiwara operators, lies midway 
between the middle and base levels. Vertices of the graph correspond to 
canonical basis elements, and Kashiwara operators remember only ``highest 
terms'' for the action of $E_i$'s on the canonical basis. This 
structure is set-theoretic. Despite its position somewhat below the first 
categorification, it is already incredibly rich. For instance, key results 
about simple and projective $R(\nu)$-modules~\cite{KL1}  follow via Kleshchev-Grojnowski 
constructions that amount to the representation-theoretical counterpart
of the crystal graph structure. 

\vspace{0.1in} 

$U^+$ is ``one-half'' of the entire quantum group $U_q(\mf{g})$. In the categorification 
of $U^+$ the diagrammatics have a braid-like behaviour, in the sense that the lines in 
diagrams only go up and do not have U-turns. Interestingly, it is possible to enlarge 
the calculus of rings $R(\nu)$ by allowing U-turns, coloring regions of resulting diagrams 
by integral weights of $\mf{g}=\mf{sl}(n)$ and adding more relations, including all 
isotopies, to categorify the entire quantum group, see Lauda~\cite{AL1} for the $\mf{sl}(2)$ case, \cite{KL3} for arbitrary $n$, and Chuang-Rouquier~\cite{CR}, Rouquier~\cite{Rou2} 
for related and less rigid axiomatics. The quantum 
group first needs to be modified, following Beilinson-Lusztig-MacPherson and 
Lusztig~\cite{Lus2},  
by adding idempotents $1_{\lambda}$ of projection onto integral weights $\lambda$. 
In this categorification generators $E_i$ and $F_i$ become biadjont functors (biadjoint 
up to grading shifts). The Grothendieck ring of the resulting 2-category is isomorphic 
to the BLM form of the quantum $\mf{sl}(n)$. The diagrammatics of biadjoint functors 
described at the beginning of this paper plays a fundamental role in the definition of
the 2-category and computations in it. Recent categorification of the $q$-Schur algebra 
by Mackaay, ~Sto\v{s}i\'c, and Vaz~\cite{MSV} 
links the Soergel category and categorified quantum $\mf{sl}(n)$.  

The $n=2$ case, due to Lauda~\cite{AL1},  settles a conjecture 
of Igor Frenkel, circa 1994, that there exists a categorification of quantum $\mf{sl}(2)$
(some motivations for the conjecture can be found in~\cite{CF}). 
Back in the early 90's Igor Frenkel, who was the author's advisor in graduate school, 
envisioned categorification lifting the entire theory of quantum groups, 
quantum 3-manifold invariants, and conformal field theory. Through the work of 
many people, his prophecy is becoming a reality.  
    
\vspace{0.1in} 
    
{\bf Acknowledgments} The author would like to thank Hiraku Nakajima and the Mathematical 
Society of Japan for the opportunity to deliver the Takagi Lectures in June 2009 and for 
their hospitality. The author is grateful to Aaron Lauda, Radmila Sazdanovi\'c, and 
Joshua Sussan for reading the manuscript and providing valuable comments and corrections. 
Partial support during the writing of present notes came from  
the NSF, via grant DMS-0706924.  The list of references is grossly inadequate and 
misses quite a few fundamental contributions - our apologies for the omissions 
go to many fine practitioners in these areas.   
    
\vspace{0.2in} 
     


\vspace{0.1in}
 
\noindent
{ \sl \small Mikhail Khovanov, Department of Mathematics, Columbia University, New York, NY 10027}

\noindent
  {\tt \small email: khovanov@math.columbia.edu}

\end{document}